\documentclass[11pt,letterpaper]{article}

\usepackage{amssymb,amsthm,amsmath}
\usepackage{graphicx}
\graphicspath{{figures/}}
\usepackage[margin=1in]{geometry}
\usepackage{color}
\usepackage{array}
\usepackage{capt-of}
\usepackage{subcaption}
\DeclareCaptionFont{introhist}{\sffamily\fontsize{7}{9}\selectfont}
\usepackage{url}
\usepackage{authblk}
\usepackage[ruled,vlined]{algorithm2e}

\newtheorem{theorem}{Theorem}[section]
\newtheorem{lemma}[theorem]{Lemma}
\newtheorem{proposition}[theorem]{Proposition}

\newtheorem{corollary}[theorem]{Corollary}

\theoremstyle{definition}

\theoremstyle{remark}

\newcommand{\bbR}{\mathbb{R}}

\newcommand{\conv}{\mathrm{conv}}

\makeatletter
\newcommand*{\rom}[1]{\expandafter\@slowromancap\romannumeral #1@}
\makeatother

\usepackage{booktabs}
\usepackage{placeins}

\title{Distributionally Robust Optimization over Incomparable Spaces}
\author[1]{Junyu Chen\thanks{Email: chenjunyu@u.nus.edu}}
\author[2]{Jieling Shi}
\author[1,3]{Yong Sheng Soh}
\affil[1]{Department of Mathematics, National University of Singapore}
\affil[2]{Department of Statistics and Data Science, National University of Singapore}
\affil[3]{Institute of Operations Research and Analytics, National University of Singapore}

\begin{document}

\maketitle

\begin{abstract}
Distributionally robust optimization (DRO) provides a principled framework for making decisions under uncertainty when the underlying distribution may be incorrectly specified. A widely used approach builds on ideas from optimal transport (OT) as a means of measuring differences between distributions. However, this depends on these distributions residing in the same space.  In certain settings, the reference data distribution we have and the space on which we make our decision may be different, which renders classical OT techniques impossible. In this paper, we provide a new DRO formulation in which ambiguity sets are defined through the Gromov–Wasserstein distance, which is an extension of OT to pairs of distributions residing in incomparable spaces. We show that several seemingly unrelated optimization formulations, including conditional value-at-risk minimization, group distributionally robust optimization, robust graph $k$-median, and $k$-median with outliers, can be recovered as special cases of our formulation. We describe the computational aspects of solving the proposed formulation based on the doubly-nonnegative relaxation. Last, we illustrate the framework through experiments on clustering over incomparable spaces and stochastic graph $k$-median. 

\end{abstract}

\section{Introduction}
\label{sec:introduction}

Stochastic programming is a modeling framework that deals with making optimal decisions under uncertainty in relevant input data.  A prominent approach is to commit to a decision that minimizes the {\em expected cost} of the associated decision whereby the expectation is taken with respect to all possible realizations of unknown parameters that are relevant to our problem.  More concretely, let's suppose that the process of making such a decision can be cast as the following stochastic program in which we seek a decision -- parameterized by $\theta \in \Theta$ -- that minimizes the following 
\begin{equation}\label{eq:decision-making}
    \min_{\theta \in \Theta} ~ \mathbb E_{x\sim\mu}[f_\theta(x)]. 
\end{equation}
Here, $x$ is a vector that models input data that is relevant to selecting a decision $\theta$ -- these can represent problem parameters, or possible outcomes of a future scenario.  In addition, $\mu$ is a probability distribution that models uncertainty in $x$, $f$ is a function that models the loss or cost of making the decision $\theta$ under a specific realization of the parameters $x$, while $\Theta$ denotes the collection of all permissible decisions $\theta$.

In many practical applications, however, one may not have access to the underlying distribution $\mu$ of the input parameters.  For instance, $\mu$ may refer to the true population distribution for which one only has access to via an empirical distribution, say $\hat{\mu}_N$, formed by a finite number of $N$ observations.   The use of an approximate or inaccurate estimate $\hat{\mu}_N$ of $\mu$ may sometimes lead to poor decisions~\cite{mohajerinesfahani2018data}. For instance, in inventory management problems, one typically sees relatively few instances of high demand during regular operations; however, understocking during a period of a high demand can be extremely costly, even if these events are rare.  To guard against these extreme scenarios, the \textit{distributionally robust optimization (DRO)} framework extends the regular stochastic programming framework as in \eqref{eq:decision-making} by selecting $\theta \in \Theta$ that minimizes the {\em worst case} risk among all possible distributions $\mu$ that belong to a pre-specified ambiguity set $\mathcal{U}$ as follows
\begin{equation}\label{eq:dro}
    \min_{\theta\in\Theta}\;
    \max_{\mu\in\mathcal{U}}\; ~
    \mathbb E_{x\sim\mu}[f_\theta(x)].
\end{equation}
The DRO formulation \eqref{eq:dro} comprises two layers.  Given a fixed $\theta \in \Theta$, the inner problem seeks the worst case $\mu \in \mathcal{U}$ that maximizes the expected cost.  The outer problem performs the optimization over $\theta \in \Theta$.  The above minimax formulation \eqref{eq:dro} has been widely studied and has found applications in inventory management \cite{scarf1958min,zhang2024optimal}, portfolio optimization \cite{delage2010distributionally,mohajerinesfahani2018data}, and supervised learning \cite{shafieezadehabadeh2019regularization,zhu2022distributionally}. Related distributionally robust formulations have also been used for hypothesis testing~\cite{gao2018robust} and power-system operation~\cite{wang2018risk}.

\subsection{Optimal Transport-based ambiguity sets}

The modeling framework \eqref{eq:dro} permits a range of ambiguity sets.  Prominent choices include those based on the moments of the distribution as well as those based on divergence measures -- an extensive list of ambiguity sets can be found in \cite[Section~2]{kuhn2025distributionally}.  In this paper we will be primarily interested in ambiguity sets that are specified via optimal transport-based distances with respect to a reference measure.  

Concretely, suppose that the input parameters $x$ are drawn from a space $\mathcal{X}$ that has metric structure specified by $d_{\mathcal{X}}$ (that is, $(\mathcal{X},d_{\mathcal{X}})$ is a metric space).  Given a pair of probability measures $\mu,\nu \in \mathcal{P}(\mathcal{X})$ over these spaces, the $p$-Wasserstein distance between these distributions is given by the minimization instance
\begin{equation} \label{eq:p-wasserstein}
	W_p (\mu,\nu) = \min_{\pi\in\Pi(\mu,\nu)} \left( \int_{\mathcal X\times\mathcal X} d_{\mathcal{X}}(x,x')^p \,\mathrm d\pi(x,x') \right)^{1/p}, \qquad p \in [1,\infty).
\end{equation}
Here, 
\begin{equation*}
	\Pi(\mu,\nu) := \{ \pi \ge 0 : \int_{x \in \mathcal X} \mathrm d\pi(x,x') = \nu, \int_{x' \in \mathcal X} \mathrm d\pi(x,x') = \mu \}
\end{equation*}
denotes the set of all feasible couplings.  
The $p$-Wasserstein-DRO is the family of ambiguity sets specified by 
\begin{equation}\label{eq:wdro-ambiguity-set}
	\mathcal{U}_\rho(\mu) = \{ \nu : \nu \in \mathcal{P}(\mathcal{X}), W_p (\mu,\nu) \le \rho \}.
\end{equation}
Under suitable assumptions, the WDRO problem admits a finite dimensional convex reformulation, and admits finite-sample performance guarantees~\cite{mohajerinesfahani2018data}.

\subsection{Optimal transport between incomparable spaces via Gromov-Wasserstein}

The formulation of the $p$-Wasserstein distance relies on access to metric evaluations $d_{\mathcal{X}}(x,x')$.  In certain settings, however, the information that one wishes to incorporate to inform the distribution $\mu\in\mathcal{P}(\mathcal{X})$ of the input parameters $x\in \mathcal{X}$ may come from a source $\mathcal{Y}$ that has similarities to $\mathcal{X}$, but is ultimately different from $\mathcal{X}$.  For instance, suppose a policymaker is deciding on a course of action during a pandemic.  The relevant input is the underlying epidemic network between the population, which informs how a disease may spread, but such information is typically difficult to obtain.  On the other hand, a good proxy for such information may available in the form of social network data. 

In such settings, the formulation of the $p$-Wasserstein distance \eqref{eq:p-wasserstein} does not ordinarily apply because one does not have access to a suitable function that describes the cost of moving mass from $x \in \mathcal{X}$ to $y \in \mathcal{Y}$ -- these spaces are fundamentally different.  The key conceptual obstruction here is: what is a principled manner to transfer knowledge that is available in the domain $\mathcal{Y}$ --  say in the form of a distribution $\nu \in \mathcal{P}(\mathcal{Y})$ -- into the domain of $\mathcal{X}$?

To this end, the Gromov-Wasserstein problem is variant of the optimal transport problem \eqref{eq:p-wasserstein} that is based on seeking couplings between pairs of distributions that are supported over two separate metric spaces so as to minimize the distortion between the underlying metric evaluations.  More formally, we call a tuple $\mathfrak{X} := (\mathcal{X},d_{\mathcal{X}},\mu_{\mathcal{X}})$ a {\em metric measure space} if the pair $(\mathcal{X},d_{\mathcal{X}})$ specifies a metric space, while $\mu_{\mathcal{X}}$ specifies a measure over $\mathcal{X}$.  Given a pair of mm-spaces $\mathfrak{X}$ and $\mathfrak{Y}$, the Gromov-Wasserstein (GW) distance between these mm-spaces is given by the solution of the following minimization instance
 \cite{memoli2011gromov}:
\begin{equation}\label{eq:gw-def}
\operatorname{GW}_{p,q} \bigl(\mathfrak{X},\mathfrak{Y}\bigr) := \inf_{\pi\in\Pi(\mu,\nu)} \Big(
 \iint_{(\mathcal{X} \times\mathcal{Y})^2}
 \bigl|d_\mathcal{X}^{q} (x,x')-d_\mathcal{Y}^{q} (y,y')\bigr|^p
 \,\mathrm d\pi(x,y)\,\mathrm d\pi(x',y') \Big)^{1/p}. 
\end{equation}
The basic intuition behind the GW formulation in \eqref{eq:gw-def} is that it compares the metric evaluations between pairs of points $(x,x')$ from the target space with metric evaluations between pairs of points $(y,y')$ from the source space, and seeks couplings that minimize the overall difference, weighted by the couplings.  In so doing, the GW provides as outputs couplings that draw comparisons between these mm-spaces without the need to specify a cost function $c$.  

Beyond these examples where a cost $c$ is absent, the GW distance is also useful in settings where it is necessary to take into account the ``local geometry'' of the problem.  A concrete illustration of such an example is in computer vision where one wishes to compute a correspondence between two point clouds representing different shapes~\cite{houry2026gromov,memoli2011gromov,solomon2016entropic}.  In this setting, the points reside in a common Euclidean space.  However, the couplings that are obtained via the formulation \eqref{eq:p-wasserstein} with $c$ being the Euclidean norm lead to transportation plans where a subset of points may ``escape the surface", which is an undesirable artifact in computer vision applications.  Couplings obtained from the GW problem typically do not suffer from such a defect as the objective explicitly penalizes for distortion in the metric space, and as a consequence obtain correspondences that largely preserve the overall metric structure. 

The GW distance and its variants have been applied to a range of tasks.  Examples of graph processing tasks include graph matching~\cite{xu2019gromov}, graph classification~\cite{vayer2019optimal} and graph representation learning~\cite{vincentcuaz2022template,vincentcuaz2021online}. Other applications include heterogeneous domain adaptation~\cite{demetci2024breaking,yan2018semi}, word alignment and machine translation~\cite{alvarezmelis2018gromov,chen2020graph}, alignment of single-cell data across modalities~\cite{cao2022pamona,demetci2022scot,ryu2025cross}, and prediction of molecular energies~\cite{peyre2016gromov}.

\subsection{Problem formulation}
\label{subsec:problem-formulation}

In this paper, we consider the task of solving stochastic programs in the form of \eqref{eq:decision-making} whereby information concerning the distribution of $x \in \mathcal{X}$ comes from a source $\mathcal{Y}$ that is different from $\mathcal{X}$.  Here, $x \in \mathcal{X}$ models the unknown parameter or uncertain outcome.  The prior knowledge comes in the form of a probability measure $\nu \in \mathcal{P}(\mathcal{Y})$ that is supported over the space $\mathcal{Y}$ that is different from $\mathcal{X}$.  Nevertheless, the information that is provided by $\nu$ remains valuable to us, and we wish to transfer knowledge from $\nu$ to our decision problem.

To address this problem, we rely on the Gromov-Wasserstein distance to describe a family of distributions that are close by on the basis of a distortion measure of the metric space.  Formally, suppose that the space $\mathcal{X}$ (describing the unknown parameters or uncertain outcomes that are relevant to our decision process) is afforded with a metric so that $(\mathcal{X},d_{\mathcal{X}})$ specifies a metric space.  We also suppose that the space where our prior information comes from $\mathcal{Y}$ also comes with metric structure so that $(\mathcal{Y},d_{\mathcal{Y}})$ specifies a metric measure space.  Suppose that the information in $\mathcal{Y}$ is presented to us in the form of a probability distribution $\nu \in \mathcal{P}(\mathcal{Y})$, which we wish to transfer to $\mathcal{X}$.  The manner in which we transfer knowledge or information from $\nu$ to our decision process \eqref{eq:decision-making} is via the DRO modeling language of \eqref{eq:dro} whereby we solve the following GW distributionally robust optimization (GWDRO) problem
\begin{equation}\label{eq:gw-dro}
    \min_{\theta\in\Theta}\;
    \max_{\mu\in\mathcal{U}_{\rho}^{\mathrm{GW}}}\;
    \mathbb{E}_{x\sim\mu}[f_{\theta}(x)],
\end{equation}
and where the ambiguity set $\mathcal{U}_{\rho}^{\mathrm{GW}}:=\mathcal{U}_{\rho}^{\mathrm{GW}}(\nu)$, is specified by a ball of measures in $\mathcal{P}(\mathcal{X})$ whose GW distance from $\nu$ is smaller than a pre-specified threshold
\begin{equation}
    \mathcal{U}_{\rho}^{\mathrm{GW}}(\nu)
    \,:=\, \left\{\mu\in \mathcal{P}(\mathcal{X}) :
    \operatorname{GW}((\mathcal{X},d_{\mathcal{X}},\mu),(\mathcal{Y},d_{\mathcal{Y}},\nu))\leq\rho\right\}.
    \label{eq:gw-ambiguity-set-cts}
\end{equation}

In particular, we envisage that the GW-based DRO formulation will be useful for decision making in the following settings:
\begin{enumerate}
\item We do not have reliable information about the distribution $\mu$ of the unknown parameters $x \in \mathcal{X}$.  Instead, our only source of information comes in the form of data over a different space $\mathcal{Y}$, which we use to form $\nu \in \mathcal{P}(\mathcal{Y})$.  The space $\mathcal{Y}$ possesses some structural similarity with $\mathcal{X}$, which suggests that information in the domain of $\mathcal{Y}$ transfers over to $\mathcal{X}$.
\item An estimate of the $\mu$ in the form of $\nu \in \mathcal{P}(\mathcal{X})$ is available.  However, it is necessary to take into consideration local geometry in the problem, which renders the ordinary WDRO framework unsuitable.
\end{enumerate}

\paragraph{GWDRO over finite metric spaces.}  In practical settings, one typically works with metric spaces that are finite -- for instance, these may be obtained as a discretization of the spaces $\mathcal{X}$ and $\mathcal{Y}$.  We specialize the above description \eqref{eq:gw-ambiguity-set-cts} to such settings.  Formally, suppose we have $\mathcal{X}=\{x_i:i\in[m]\}$ as the set of all possible parameters or outcomes, and we have $\mathcal{Y}=\{y_k:k\in[n]\}$.  We let $\Delta_k := \{ p \in \mathbb{R}_{+}^{k} : \sum_{i\in[k]}p_i=1 \}$ denote the probability simplex. Then the set of probability distributions $\mathcal{P}(\mathcal{X})$ and $\mathcal{P}(\mathcal{Y})$ can be identified with $\Delta_m$ and $\Delta_n$, respectively.  Let $\mu_{\mathcal{X}}$ and $\mu_{\mathcal{Y}}$ denote the unknown target and source distributions.  The set of all feasible couplings is expressed as set of finite dimensional matrices $\{ \pi \in \mathbb{R}^{m \times n}_{+} : \pi\mathbf{1}_n=\mu_{\mathcal{X}}, \pi^{\top}\mathbf{1}_m = \mu_{\mathcal{Y}} \}$, where $\mathbf{1}_r$ denotes the vector of ones in $\mathbb{R}^{r}$.  The metric spaces $(\mathcal{X},d_{\mathcal{X}})$ and $(\mathcal{Y},d_{\mathcal{Y}})$ can be conveyed via finite dimensional matrices $D^{\mathcal{X}} \in \mathbb{R}^{m \times m}$ where $D^{\mathcal{X}}_{ij} := d^q_{\mathcal{X}}(x_i,x_j)$, $i,j\in[m]$, and $D^{\mathcal{Y}} \in \mathbb{R}^{n \times n}$, where $D^{\mathcal{Y}}_{kl} := d^q_{\mathcal{Y}}(y_k,y_l)$, $k,l\in[n]$.  The GW objective is subsequently given by the following quadratic function
\begin{equation}
    \mathcal{D}(\pi)
    :=
    \sum_{i,j\in[m]}\sum_{k,l\in[n]}
    \bigl(D^{\mathcal{X}}_{ij}-D^{\mathcal{Y}}_{kl}\bigr)^2
    \pi_{ik}\pi_{jl}.
    \label{eq:gw-distortion}
\end{equation}
Suppose we make the choice $p=2$ within \eqref{eq:gw-def}.  Then the GW ambiguity set is given by the set of all probability vectors that satisfy the following condition
\begin{equation}
\mathcal{U}_{\rho}^{\mathrm{GW}}(\widehat{\mu}_{\mathcal{Y}})
=\left\{\mu\in\Delta_m:\ \exists\,\pi\in\mathbb{R}_{+}^{m\times n},\ \pi\mathbf{1}_n=\mu,\ \pi^{\top}\mathbf{1}_m=\widehat{\mu}_{\mathcal{Y}},\ \mathcal{D}(\pi)\leq\rho^2\right\}.
\label{eq:finite-gw-ambiguity-set}
\end{equation} 

The target distribution is determined by the row sums of $\pi$.
Thus, for every feasible coupling,
\begin{equation*}
    \mathbb{E}_{x\sim\mu}[f_{\theta}(x)]
    =
    \sum_{i\in[m]}\mu_i f_{\theta}(x_i)
    =
    \sum_{i\in[m]}\sum_{k\in[n]}f_{\theta}(x_i)\pi_{ik}.
\end{equation*}
We can therefore eliminate $\mu$ and rewrite \eqref{eq:gw-dro} as
\begin{equation}\tag{GW-DRO}
\label{eq:special-cases-gwdro}
\min_{\theta\in\Theta}\;
\max_{\substack{
\pi\in\mathbb{R}_{+}^{m\times n}\\
\mathcal{D}(\pi)\leq\rho^2\\
\pi^{\top}\mathbf{1}_m=\widehat{\mu}_{\mathcal{Y}}}}
\sum_{i\in[m]}\sum_{k\in[n]}f_{\theta}(x_i)\pi_{ik}. 
\end{equation}

\paragraph{Outline.}  In Section~\ref{sec:gwdro-special-cases}, we discuss a series of seemingly unrelated statistical problems (CVaR, $k$-medians, group DRO, semi-GW) that can in fact be interpreted as special instances of the GW-DRO formulation \eqref{eq:gw-dro}.  In Section~\ref{sec:worst-case-upper-bound} we discuss computational aspects of solving \eqref{eq:gw-dro}.  In Section~\ref{sec:numerical-experiments} we implement our framework on numerical experiments involving synthetic data.

\section{Examples}
\label{sec:gwdro-special-cases}

In this section we discuss a list of seemingly unrelated statistical problems that can in fact be understood as special instances of our proposed framework \eqref{eq:special-cases-gwdro}.  These connections point towards the modeling capabilities of \eqref{eq:special-cases-gwdro}.  To keep the discussion brief, we only state the main propositions that formally connect the existing statistical models to our framework -- the accompanying proofs are provided in Appendix~\ref{app:proofs-examples}.

\subsection{Semi-relaxed GW}
\label{subsec:special-cases-srgw}

This subsection is organized as follows. First, we introduce semi-relaxed GW (srGW) and describe how it can be used as a benchmark for decision making over incomparable spaces, where the source distribution that provides the available information is defined on a space incomparable with the target space. Roughly speaking, the srGW problem estimates a target distribution by projecting the source distribution onto the target space under the GW distance. One can then make decisions under the estimated target distribution. Second, we show that, when the srGW estimate is unique, the resulting decision-making problem can be recovered from \eqref{eq:special-cases-gwdro} by choosing the smallest radius for which the GW ambiguity set is nonempty. Last, we consider WDRO with a Wasserstein ambiguity set centered at the srGW estimate and use an example to compare this ambiguity set with the GW ambiguity set. This motivates the benefit of using the GW ambiguity set. 

We first introduce the srGW problem and describe how it can be used to make decisions under an estimated target distribution. The srGW problem considers the task of choosing a target distribution to minimize its GW distance from a given source distribution \cite{vincentcuaz2022semi}. More precisely, let $(\mathcal{X},d_{\mathcal{X}})$ and $(\mathcal{Y},d_{\mathcal{Y}})$ be the finite metric spaces introduced in Section~\ref{subsec:problem-formulation}, with $|\mathcal{X}|=m$ and $|\mathcal{Y}|=n$. For a source distribution $\nu\in\Delta_n$, the srGW problem is
\begin{equation}
\label{eq:special-cases-srgw-radius}
\rho_{\min}(\nu)^2:=\min_{\mu\in\Delta_m}\operatorname{GW}_{2}^{2}(\mu,\nu).
\end{equation}
The srGW problem can be viewed as the projection of the source metric measure space $(\mathcal{Y},d_{\mathcal{Y}},\nu)$ onto the family $\{(\mathcal{X},d_{\mathcal{X}},\mu):\mu\in\Delta_m\}$ with respect to the GW distance~\cite{vincentcuaz2022semi}. The srGW formulation has been studied for graph partitioning, graph clustering, dictionary learning, and graph completion~\cite{vincentcuaz2022semi}. 

In our setting, srGW can be used to estimate target probabilities from the observed source distribution. For the observed source distribution $\widehat{\mu}_{\mathcal{Y}}$, let $\widehat{\mu}_{\mathcal{X}}^{\mathrm{srGW}}$ be an optimizer of srGW, that is
\begin{equation*}
\widehat{\mu}_{\mathcal{X}}^{\mathrm{srGW}}\in\operatorname*{argmin}_{\mu\in\Delta_m}\operatorname{GW}_{2}^{2}\bigl(\mu,\widehat{\mu}_{\mathcal{Y}}\bigr).
\end{equation*}
To choose a decision that minimizes the expected loss under $\widehat{\mu}_{\mathcal{X}}^{\mathrm{srGW}}$, it suffices to solve
\begin{equation}
\label{eq:special-cases-srgw-formulation}
\min_{\theta\in\Theta}\sum_{i\in[m]}\widehat{\mu}_{\mathcal{X},i}^{\mathrm{srGW}}f_{\theta}(x_i).
\end{equation}

We next show that, when the minimizer $\widehat{\mu}_{\mathcal{X}}^{\mathrm{srGW}}$ is unique, problem \eqref{eq:special-cases-srgw-formulation} can be recovered from \eqref{eq:special-cases-gwdro} by choosing the smallest radius for which the GW ambiguity set centered at $\widehat{\mu}_{\mathcal{Y}}$ is nonempty.

\begin{proposition}
\label{prop:special-cases-srgw}
Suppose that \eqref{eq:special-cases-srgw-radius} with $\nu=\widehat{\mu}_{\mathcal{Y}}$ has a unique minimizer $\widehat{\mu}_{\mathcal{X}}^{\mathrm{srGW}}$. Then problem \eqref{eq:special-cases-gwdro} at radius $\rho=\rho_{\min}(\widehat{\mu}_{\mathcal{Y}})$ is equivalent to problem \eqref{eq:special-cases-srgw-formulation}.
\end{proposition}

Proposition~\ref{prop:special-cases-srgw} establishes the equivalence at the smallest feasible radius when the minimizing target distribution is unique. Increasing the radius in \eqref{eq:special-cases-gwdro} allows additional adversarial target distributions. We next consider an alternative way to account for uncertainty: we construct a Wasserstein ambiguity set centered at the srGW estimate $\widehat{\mu}_{\mathcal{X}}^{\mathrm{srGW}}$ and compare it with the GW ambiguity set. The resulting WDRO problem is
\begin{equation}\label{eq:wdro-incomparable}
    \min_{\theta\in\Theta}\;
    \max_{\mu\in\mathcal{U}^{\mathrm{srGW}}_\rho}\;
    \mathbb{E}_{{x}\sim\mu}[f_{\theta}({x})],
\end{equation}
where
\begin{equation}\label{eq:srGW-ambiguity}
\mathcal{U}^{\mathrm{srGW}}_\rho = \{\nu \in \Delta_m: W_2^2(\nu,\widehat{\mu}_{\mathcal{X}}^{\mathrm{srGW}}) \le \rho^2 \}
\end{equation}
is the Wasserstein ambiguity set with center $\widehat{\mu}_{\mathcal{X}}^{\mathrm{srGW}}$.

\begin{figure}[!htbp]
    \centering
    \captionsetup[subfigure]{font=introhist,justification=centering,singlelinecheck=false,skip=2pt}
    \begin{subfigure}[t]{0.195\linewidth}
        \centering
        \includegraphics[width=\linewidth]{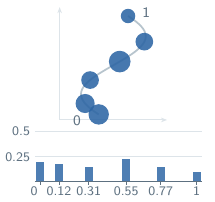}
        \caption{Observed source\\$\widehat{\mu}_{\mathcal{Y},N}$}
        \label{fig:intro-hist-observed-source}
    \end{subfigure}\hfill
    \begin{subfigure}[t]{0.195\linewidth}
        \centering
        \includegraphics[width=\linewidth]{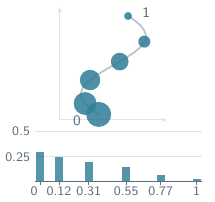}
        \caption{True source\\$\mu_{\mathcal{Y}}$}
        \label{fig:intro-hist-true-source}
    \end{subfigure}\hfill
    \begin{subfigure}[t]{0.195\linewidth}
        \centering
        \includegraphics[width=\linewidth]{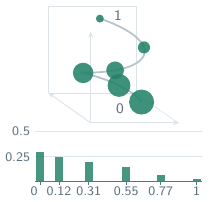}
        \caption{True target\\$\mu_{\mathcal{X}}$}
        \label{fig:intro-hist-true-target}
    \end{subfigure}\hfill
    \begin{subfigure}[t]{0.195\linewidth}
        \centering
        \includegraphics[width=\linewidth]{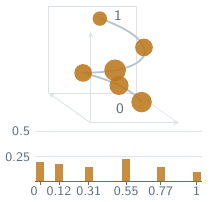}
        \caption{srGW estimate\\$\widehat{\mu}_{\mathcal{X}}^{\mathrm{srGW}}$}
        \label{fig:intro-hist-srgw-estimate}
    \end{subfigure}\hfill
    \begin{subfigure}[t]{0.195\linewidth}
        \centering
        \includegraphics[width=\linewidth]{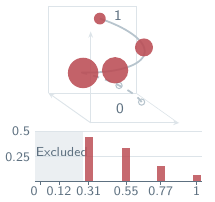}
        \caption{Conditional distribution\\$\mu_{\mathcal{X}}^{\mathrm{cond}}$}
        \label{fig:intro-hist-conditional-target}
    \end{subfigure}
    \begingroup
    \small
    \caption[Source sampling error and target ambiguity sets]{Distributions on source and target curves.
    Dot areas and bar heights represent probabilities; histogram axes show normalized arc positions. The dashed curve and shaded region in \subref{fig:intro-hist-conditional-target} indicate the segment excluded by conditioning. The Wasserstein ambiguity set centered at $\widehat{\mu}_{\mathcal{X}}^{\mathrm{srGW}}$ cannot contain $\mu_{\mathcal{X}}$ while excluding $\mu^{\mathrm{cond}}_{\mathcal{X}}$, whereas the GW ambiguity set centered at $\widehat{\mu}_{\mathcal{Y},N}$ can.}
    \label{fig:intro-gw-wasserstein-histograms}
    \endgroup
\end{figure}

We show that the Wasserstein ambiguity set centered at the srGW estimate, however, has some shortcomings for decision making over incomparable spaces. Specifically, if the empirical source distribution is a poor estimate of the underlying source distribution, the srGW estimate of the target distribution may also be inaccurate. The following example shows that increasing the Wasserstein radius to include the true target distribution can also require including a distribution that gives a misleading description of the target distribution. See Figure \ref{fig:intro-gw-wasserstein-histograms} for an illustration. In this example, the underlying source distribution $\mu_{\mathcal{Y}}$ in \ref{fig:intro-gw-wasserstein-histograms}\subref{fig:intro-hist-true-source} is supported on six points on a curve. What the decision maker observes is the empirical source distribution $\widehat{\mu}_{\mathcal{Y},N}$ in \ref{fig:intro-gw-wasserstein-histograms}\subref{fig:intro-hist-observed-source}, obtained by sampling $N=40$ points according to the true source distribution. The decision maker wants to make decisions under uncertainty on the target space in \ref{fig:intro-gw-wasserstein-histograms}\subref{fig:intro-hist-true-target}, where the possible outcomes are known but their probabilities are unknown. In fact, the underlying target distribution $\mu_{\mathcal{X}}$ in \ref{fig:intro-gw-wasserstein-histograms}\subref{fig:intro-hist-true-target} is obtained by bending the true source distribution in \ref{fig:intro-gw-wasserstein-histograms}\subref{fig:intro-hist-true-source} and embedding it in $\bbR^3$. Following the framework in \eqref{eq:wdro-incomparable}, the decision maker transfers the observed source distribution in \ref{fig:intro-gw-wasserstein-histograms}\subref{fig:intro-hist-observed-source} to the target space using srGW and obtains $\widehat{\mu}_{\mathcal{X}}^{\mathrm{srGW}}$, as illustrated in \ref{fig:intro-gw-wasserstein-histograms}\subref{fig:intro-hist-srgw-estimate}. In addition, consider a constructed undesired distribution $\mu^{\mathrm{cond}}_\mathcal{X}$, which is obtained by restricting the true target distribution to the segment $[0.31,1]$ when sampling. Suppose that the decision maker constructs the Wasserstein uncertainty set $\mathcal{U}^{\mathrm{srGW}}_\rho$ as defined in \eqref{eq:srGW-ambiguity}. Note that $W_2^2(\mu_\mathcal{X},\widehat{\mu}_{\mathcal{X}}^{\mathrm{srGW}}) = 0.03396 > W_2^2(\mu^{\mathrm{cond}}_\mathcal{X},\widehat{\mu}_{\mathcal{X}}^{\mathrm{srGW}}) = 0.03329$. Therefore, if $\rho$ is chosen to be sufficiently large ($\rho^2 \ge 0.03396$) for $\mathcal{U}^{\mathrm{srGW}}_\rho$ to include the underlying target distribution, then $\mathcal{U}^{\mathrm{srGW}}_\rho$ will also include the undesired distribution. Conversely, if $\rho$ is chosen to be small enough ($\rho^2 < 0.03329$) for $\mathcal{U}^{\mathrm{srGW}}_\rho$ to exclude the undesired distribution, then $\mathcal{U}^{\mathrm{srGW}}_\rho$ will also exclude the real target distribution. Thus, no choice of the Wasserstein radius includes the true target distribution in \ref{fig:intro-gw-wasserstein-histograms}\subref{fig:intro-hist-true-target} while excluding the conditional distribution in \ref{fig:intro-gw-wasserstein-histograms}\subref{fig:intro-hist-conditional-target}.

In contrast, consider the GW ambiguity set in \eqref{eq:finite-gw-ambiguity-set} with source distribution $\widehat{\mu}_{\mathcal{Y},N}$. The squared GW distances satisfy $\operatorname{GW}_2^2(\mu_{\mathcal{X}},\widehat{\mu}_{\mathcal{Y},N})\leq0.024$ and $\operatorname{GW}_2^2(\mu^{\mathrm{cond}}_{\mathcal{X}},\widehat{\mu}_{\mathcal{Y},N})\geq0.035$. Thus, for $0.024\leq\rho^2<0.035$, the GW ambiguity set contains the true target distribution and excludes the conditional distribution. 

The difference arises from how the two ambiguity sets use the relation between the source and target spaces. In the srGW-based approach, this relation is used to construct a single estimate of the target distribution. Once the estimate is obtained, the Wasserstein constraint measures the cost of transporting probability between the estimate and each candidate target distribution. The constraint does not directly measure how well the candidate distribution preserves the within-space distances of the empirical source distribution under a coupling. Consequently, a candidate distribution can be closer to the srGW transferred estimate in Wasserstein distance while having a larger GW distance from the empirical source distribution. Increasing the Wasserstein radius does not change this ordering. In contrast, the GW constraint evaluates every candidate target distribution through its coupling with the empirical source distribution and the associated distortion of within-space distances. When source observations inform target probabilities through a structural relationship between the two spaces, this comparison allows that relationship to be imposed directly in the ambiguity set. The preceding example illustrates the resulting distinction: a GW ambiguity set can include the true target distribution while excluding the conditional distribution, whereas no Wasserstein ball centered at the srGW estimate can do so. This motivates constructing the ambiguity set using the GW distance.

\subsection{Conditional value-at-risk}
\label{subsec:special-cases-cvar}

The conditional value-at-risk (CVaR) is a risk measure defined as the expected loss in the upper tail of a distribution.  In simpler words, it is the average loss conditioned on the event that a loss occurred.  Given a random variable $Z$ which models loss as well as a confidence level $\alpha\in(0,1)$, the CVaR is given by the solution of the following \cite[Theorem~10]{rockafellar2002conditional}
\begin{equation*}
\operatorname{CVaR}_{\alpha}(Z)
=
\min_{t\in\mathbb{R}}
\left\{
t+\frac{1}{1-\alpha}\mathbb{E}\bigl[(Z-t)_{+}\bigr]
\right\}.
\end{equation*}
Here $(a)_{+}=\max\{a,0\}$.  In practice, the CVaR is often applied to $Z=f_{\theta}(x)$, where $f_{\theta}$ is a loss function with parameter $\theta\in\Theta$ and $x\sim\operatorname{Unif}(x_1,\ldots,x_m)$ is an uncertain outcome taking values uniformly form $\mathcal{X} = \{x_1,\ldots,x_m\}$~\cite{rockafellar2000optimization}. In this case, minimizing CVaR takes the form
\begin{equation}
\label{eq:special-cases-cvar-formulation}
\min_{\theta\in\Theta}\operatorname{CVaR}_{\alpha}\bigl(f_{\theta}(x)\bigr)
=
\min_{\theta\in\Theta,\,t\in\mathbb{R}}
\left\{
t+\frac{1}{(1-\alpha)m}\sum_{i\in[m]}\bigl(f_{\theta}(x_i)-t\bigr)_{+}
\right\}.
\end{equation}

In what follows, we prove that \eqref{eq:special-cases-cvar-formulation} with $\alpha=1-n/m$ can be recovered from \eqref{eq:special-cases-gwdro}. To this end, we first write \eqref{eq:special-cases-cvar-formulation} as a worst-case problem over an ambiguity set and describe this set through its extreme points. We then choose source and target metric spaces for which the GW ambiguity set consists of these extreme points. 

By applying the dual representation of CVaR in~\cite[Chapter~6]{shapiro2014lectures} to \eqref{eq:special-cases-cvar-formulation}, we have the following.

\begin{lemma}
\label{lem:special-cases-cvar-dual}
Let $n\in[m-1]$, $\alpha=1-n/m$, and define
\begin{equation*}
\mathcal{C}_n
:=
\bigl\{\mu\in\Delta_m:\ \mu_i\leq 1/n\ \text{for all}\ i\in[m]\bigr\}.
\end{equation*}
Then, for every $\theta\in\Theta$, one has that
\begin{equation*}
\operatorname{CVaR}_{\alpha}\bigl(f_{\theta}(x)\bigr)
=
\max_{\mu\in\mathcal{C}_n}\sum_{i\in[m]}\mu_i f_{\theta}(x_i).
\end{equation*}
\end{lemma}

Since the objective $\sum_{i\in[m]}\mu_i f_{\theta}(x_i)$ is linear in $\mu$, its maximum over $\mathcal{C}_n$ is attained at an extreme point of $\mathcal{C}_n$. We therefore describe these extreme points.  Given a subset $S\subseteq\mathcal{X}$, let $\operatorname{Unif}(S)$ denote the distribution that assigns probability $1/|S|$ to each point in $S$ and zero probability to each point outside $S$. The scaled set $n\,\mathcal{C}_n=\{x\in[0,1]^m:\ \sum_{i\in[m]}x_i=n\}$ is the hypersimplex, which is the convex hull of the $0/1$ vectors with exactly $n$ entries equal to one~\cite{lam2007alcoved}. Dividing by $n$ gives the following lemma.

\begin{lemma}
\label{lem:special-cases-capped-simplex-hull}
For every integer $n\in[m]$,
\begin{equation*}
\mathcal{C}_n=\conv\bigl\{\operatorname{Unif}(S):\ S\subseteq\mathcal{X},\ |S|=n\bigr\}.
\end{equation*}
\end{lemma}

By Lemmas~\ref{lem:special-cases-cvar-dual} and~\ref{lem:special-cases-capped-simplex-hull}, and since the maximum of a linear function over a convex hull equals its maximum over the generating set, it suffices to choose source and target metric spaces whose GW ambiguity set consists of the distributions $\operatorname{Unif}(S)$ with $|S|=n$.

\begin{proposition}
\label{prop:special-cases-cvar}
Let $f_{\theta}:\mathcal{X}\to\mathbb{R}$ be the cost function for the decision $\theta\in\Theta$. Suppose $2\leq n<m$, $\widehat{\mu}_{\mathcal{Y}}=\mathbf{1}_n/n$, and
\begin{equation*}
D^{\mathcal{Y}}_{kl}
=
\begin{cases}
0,&k=l,\\
1,&k\neq l,
\end{cases}
\qquad
D^{\mathcal{X}}_{ij}
=
\begin{cases}
0,&i=j,\\
1,&i\neq j.
\end{cases}
\end{equation*}
Set $\rho=0$ in \eqref{eq:special-cases-gwdro}.  Then problem \eqref{eq:special-cases-gwdro} is equivalent to problem \eqref{eq:special-cases-cvar-formulation} with $\alpha=1-n/m$.
\end{proposition}

\subsection{Group distributionally robust optimization}
\label{subsec:special-cases-group-dro}

The group distributionally robust optimization problem considers the task of choosing a decision to minimize the largest expected loss over a given collection of groups~\cite{sagawa2020distributionally}. More precisely, let $\mathcal{X}=\{x_i:i\in[m]\}$ be partitioned into $G\geq1$ given, nonempty groups $\mathcal{X}_g$, $g\in[G]$. For a decision $\theta\in\Theta$, let $f_{\theta}(x)$ be the loss at $x$, as defined in Section~\ref{subsec:problem-formulation}. The empirical average loss in group $\mathcal{X}_g$ is $|\mathcal{X}_g|^{-1}\sum_{x\in\mathcal{X}_g}f_{\theta}(x)$. The empirical group distributionally robust optimization problem can be formulated as
\begin{equation}
\label{eq:special-cases-group-dro-formulation}
\min_{\theta\in\Theta}\;\max_{g\in[G]}\frac{1}{|\mathcal{X}_g|}\sum_{x\in\mathcal{X}_g}f_{\theta}(x).
\end{equation}

In what follows, we prove that, when $|\mathcal{X}_g|=n$ for every $g\in[G]$, problem~\eqref{eq:special-cases-group-dro-formulation} can be recovered from \eqref{eq:special-cases-gwdro}. As in Section~\ref{subsec:special-cases-cvar}, we take $\rho=0$, so the feasible target distributions are copies of the uniform source distribution on $n$ distinct target points that preserve all pairwise distances. We choose the target distances so that $n$ points are pairwise at the source distance $1$ only when they lie in a single group.

\begin{proposition}
\label{prop:special-cases-group-dro}
Suppose $|\mathcal{X}_g|=n\geq2$ for every $g\in[G]$, so that $m=Gn$. Let the source space be $\mathcal{Y}=\{y_k:k\in[n]\}$, with probabilities $\widehat{\mu}_{\mathcal{Y}}=\mathbf{1}_n/n$ and metric $d_{\mathcal{Y}}(y_k,y_l)=1$ for $k\neq l$. Equip the target space $\mathcal{X}$ with the metric
\begin{equation*}
d_{\mathcal{X}}(x,x')
=
\begin{cases}
0,&x=x',\\
1,&x\neq x'\text{ and }x,x'\in\mathcal{X}_g
       \text{ for some }g,\\
2,&x\in\mathcal{X}_g,\ x'\in\mathcal{X}_{h},\ g\neq h.
\end{cases}
\end{equation*}
Set $\rho=0$ in \eqref{eq:special-cases-gwdro}.  Then problem \eqref{eq:special-cases-gwdro} is equivalent to problem \eqref{eq:special-cases-group-dro-formulation}.
\end{proposition}

\subsection{Graph $k$-median}
\label{subsec:special-cases-graph-k-median}

The graph $k$-median problem considers the task of selecting $k$ centers on a graph to minimize the sum of distances from each vertex to their nearest centers, weighted by the weight of each vertex. More precisely, let $V$ be the vertex set of a finite connected undirected graph with positive edge lengths, and let $d$ be the shortest-path metric on $V$. Fix an integer $1\leq k\leq|V|$. For a set of centers $C\subseteq V$ with $|C|=k$, define the distance from $v\in V$ to its nearest center as
\begin{equation}
\label{eq:special-cases-outlier-center-distance}
d(v,C):=\min_{c\in C}d(v,c).
\end{equation}
Suppose that a vertex $v\in V$ is drawn with probability $p_v$, where $p_v\geq0$ and $\sum_{v\in V}p_v=1$. The stochastic graph $k$-median problem can be formulated as:
\begin{equation}
\label{eq:special-cases-stochastic-k-median}
\min_{\substack{C\subseteq V\\|C|=k}}\sum_{v\in V}p_v d(v,C).
\end{equation}
We consider two robust variants of the graph $k$-median problem and show that both can be recovered from \eqref{eq:special-cases-gwdro}.

\paragraph{Robust $k$-median.}

Let $S_1,\ldots,S_G$ be distinct subsets of $V$, where $G\geq1$. The robust $k$-median problem~\cite{anthony2010plant} chooses $k$ centers to minimize the largest cost over the sets $S_1,\ldots,S_G$:
\begin{equation}
\label{eq:special-cases-robust-k-median}
\min_{\substack{C\subseteq V\\|C|=k}}\;\max_{i\in[G]}\sum_{v\in S_i}d(v,C).
\end{equation}
In particular, when $S_1,\ldots,S_G$ are the singleton subsets of $V$, the cost for a given $C$ is $\max_{v\in V}d(v,C)$. Thus, in this case, the robust k-median problem replaces the expected distance in \eqref{eq:special-cases-stochastic-k-median} with the largest distance.

In what follows, we prove that problem \eqref{eq:special-cases-robust-k-median} can be recovered from \eqref{eq:special-cases-gwdro} by taking each set $S_i$ to be a singleton and choosing a sufficiently large radius.

\begin{proposition}
\label{prop:special-cases-robust-k-median}
Let $\mathcal{X}=\{S_i:i\in[G]\}$ be equipped with any metric $d_{\mathcal{X}}$, and let $(\mathcal{Y},d_{\mathcal{Y}})$ be any finite source metric space with $\mathcal{Y}=\{y_l:l\in[n]\}$ and empirical probabilities $\widehat{\mu}_{\mathcal{Y}}\in\Delta_n$. Take $\Theta=\{C\subseteq V:|C|=k\}$ and $f_{\theta}(S_i)=\sum_{v\in S_i}d(v,\theta)$ for every $\theta\in\Theta$ and $i\in[G]$. Define
\begin{equation}
\label{eq:special-cases-source-radius}
\rho(\widehat{\mu}_{\mathcal{Y}})
:=\left(\sum_{k,l\in[n]}\widehat{\mu}_{\mathcal{Y},k}\widehat{\mu}_{\mathcal{Y},l}\bigl(D^{\mathcal{Y}}_{kl}\bigr)^2\right)^{1/2}.
\end{equation}
For every $\rho\geq\rho(\widehat{\mu}_{\mathcal{Y}})$, problem \eqref{eq:special-cases-gwdro} is equivalent to problem \eqref{eq:special-cases-robust-k-median}.
\end{proposition}

We discuss the relation between \eqref{eq:special-cases-gwdro} adpated for the stochastic graph $k$-median problem in \eqref{eq:special-cases-stochastic-k-median} and the robust graph $k$-median problem in \eqref{eq:special-cases-robust-k-median} using Proposition \ref{prop:special-cases-robust-k-median}. Let $V=\{v_i:i\in[|V|]\}$, and consider the case where $S_i=\{v_i\}$ for every $i\in[|V|]$. In this case, the robust graph $k$-median problem in \eqref{eq:special-cases-robust-k-median} becomes
\begin{equation}
\label{eq:special-cases-singleton-robust-k-median}
\min_{\substack{C\subseteq V\\|C|=k}}\;\max_{v\in V}d(v,C).
\end{equation}
Set the target and source metric spaces to be $(\mathcal{X},d_\mathcal{X}) = (\mathcal{Y},d_\mathcal{Y}) = (V,d)$, where $d$ is the shortest path distance for the graph, and let $p$ be the source distribution. Then \eqref{eq:special-cases-gwdro} reads
\begin{equation}
\label{eq:special-cases-singleton-gwdro-k-median}
\min_{\substack{C\subseteq V\\|C|=k}}\;
\max_{\mu \in\mathcal{U}_{\rho}^{\mathrm{GW}}(p)}
\sum_{v\in V}\mu_v d(v,C),
\end{equation}
where, by the definition of the GW ambiguity set in \eqref{eq:finite-gw-ambiguity-set}, the GW ambigity set is given by
\begin{equation*}
\mathcal{U}_{\rho}^{\mathrm{GW}}(\widehat{\mu}_{\mathcal{Y}})
=\left\{\pi\mathbf{1}_n:\pi\in\mathbb{R}_{+}^{|V|\times |V|},\,\pi^{\top}\mathbf{1}_{|V|}=p,\,\sum_{\substack{i,j\in[|V|]\\k,l\in[|V|]}}\bigl(d(v_i,v_j)-d(v_k,v_l)\bigr)^2\pi_{ik}\pi_{jl}\leq\rho^2\right\}.
\end{equation*}
Thus, \eqref{eq:special-cases-singleton-gwdro-k-median} can be viewed as a distributionally robust version of the stochastic graph $k$-median problem in \eqref{eq:special-cases-stochastic-k-median}, obtained by considering the adversarial distributions in the GW ambiguity set centered at the source distribution $p$.

Proposition~\ref{prop:special-cases-robust-k-median} shows that, if $\rho\geq\rho(\widehat{\mu}_{\mathcal{Y}})$, with $\rho(\widehat{\mu}_{\mathcal{Y}})$ defined in \eqref{eq:special-cases-source-radius}, problem \eqref{eq:special-cases-singleton-gwdro-k-median} is equivalent to problem \eqref{eq:special-cases-singleton-robust-k-median}. In Section~\ref{subsec:branch-permutation-k-median}, we compare the optimal value of \eqref{eq:special-cases-singleton-gwdro-k-median} with the optimal value of \eqref{eq:special-cases-singleton-robust-k-median} as $\rho$ varies, and examine the tightness of the radius $\rho(\widehat{\mu}_{\mathcal{Y}})$ for the instance considered there.

\paragraph{$k$-median with outliers.}

Let $z$ be an integer satisfying $0\leq z<|V|$. The $k$-median problem with outliers chooses $k$ centers and discards $z$ vertices to minimize the total distance from the remaining vertices to their nearest centers~\cite[Section~7.3]{gupta2025structural}. It can be formulated as the following optimization problem:
\begin{equation}
\label{eq:special-cases-k-median-outliers}
\min_{\substack{C\subseteq V\\|C|=k}}\;\min_{\substack{O\subseteq V\\|O|= z}}\sum_{v\in V\setminus O}d(v,C).
\end{equation}
In particular, when the number of vertices to be discarded $z=0$, objective function of \eqref{eq:special-cases-k-median-outliers} reduces to the stochastic objective in \eqref{eq:special-cases-stochastic-k-median} with $p_v=1/|V|$.

We first apply strong duality for linear programming to rewrite the minimization over $O$ in \eqref{eq:special-cases-k-median-outliers} as a maximization over its dual variable $t \in \bbR$. Then we construct the metric spaces $(\mathcal{X},d_\mathcal{X})$ and $(\mathcal{Y},d_\mathcal{Y})$, and select the radius $\rho$ of the GW ambiguity set such that the worst-case problem in \eqref{eq:special-cases-gwdro} is equivalent to the dual program. 

\begin{lemma}
\label{lem:special-cases-outlier-duality}
For every nonempty $C\subseteq V$ and integer $0\leq z<|V|$, it holds that
\begin{equation}
\label{eq:special-cases-outlier-duality}
\min_{\substack{O\subseteq V\\|O|= z}}\sum_{v\in V\setminus O}d(v,C)
=\max_{t\in\mathbb{R}}\left\{(|V|-z)t-\sum_{v\in V}\bigl(t-d(v,C)\bigr)_{+}\right\}.
\end{equation}
Moreover, the maximization admits an optimal solution $t^{*}\in\{d(v,C):v\in V\}$.
\end{lemma}

Motivated by Lemma~\ref{lem:special-cases-outlier-duality}, we define the cost function
\begin{equation}
\label{eq:special-cases-outlier-threshold-cost}
f_C(t):=(|V|-z)t-\sum_{v\in V}\bigl(t-d(v,C)\bigr)_{+},\qquad t\in\mathbb{R}.
\end{equation}
Next, we specify the target metric space $(\mathcal{X},d_{\mathcal{X}})$ and the source metric space $(\mathcal{Y},d_{\mathcal{Y}})$ used to establish the equivalence of \eqref{eq:special-cases-k-median-outliers} and \eqref{eq:special-cases-gwdro}. We take the target points to be values of the variable $t$ at which $f_C$ can attain its maximum. Since $C\subseteq V$, every distance $d(v,C)$ belongs to the finite set
\begin{equation}
\label{eq:special-cases-outlier-threshold-space}
\mathcal{X}:=\{d(v,w):v,w\in V\}.
\end{equation}
Set $m:=|\mathcal{X}|$ and denote the distinct elements of $\mathcal{X}$ by $t_1,\ldots,t_m$.
We equip $\mathcal{X}$ with the metric $d_{\mathcal{X}}(t,t')=|t-t'|$. The source space $(\mathcal{Y},d_{\mathcal{Y}})$ may be any finite metric space.

\begin{proposition}
\label{prop:special-cases-k-median-outliers}
Let $\mathcal{X}$ be defined as \eqref{eq:special-cases-outlier-threshold-space}, equipped with $d_{\mathcal{X}}(t,t')=|t-t'|$, and let $(\mathcal{Y},d_{\mathcal{Y}})$ be any finite source metric space with $|\mathcal{Y}|=n$ and empirical probabilities $\widehat{\mu}_{\mathcal{Y}}\in\Delta_n$. Take $\Theta=\{C\subseteq V:|C|=k\}$ and the cost $f_C$ defined in \eqref{eq:special-cases-outlier-threshold-cost}. For every $\rho\geq\rho(\widehat{\mu}_{\mathcal{Y}})$, with $\rho(\widehat{\mu}_{\mathcal{Y}})$ defined in \eqref{eq:special-cases-source-radius}, problem \eqref{eq:special-cases-gwdro} with cost $f_C$ is equivalent to problem \eqref{eq:special-cases-k-median-outliers}.
\end{proposition}

\section{A procedure for computing \eqref{eq:special-cases-gwdro} via DNN relaxations}
\label{sec:worst-case-upper-bound}

In this section we discuss computational aspects pertaining to solving  \eqref{eq:special-cases-gwdro}.  There are two conceptual difficulties that arise within the formulation \eqref{eq:special-cases-gwdro}, as it is stated.  First, the inner worst-case optimization instance is stated in terms of a Gromov-Wasserstein (GW) ball.  This is problematic because the GW distance between a pair of mm-spaces is expressed as the solution of a non-convex quadratically constrained quadratic program (QCQP), and is NP-Hard to compute in general  \cite{BGP:26}.  Second, the outer minimization requires optimizing for an objective that is expressed indirectly as the solution of a maximization problem, which also introduces additional challenges.

To address the first difficulty, we describe a doubly non-negative (DNN) relaxation of the inner worst-case problem.  As a reminder, a matrix $X \in \mathbb{S}^{n}$ is said to be {\em doubly non-negative (DNN)} if it is positive semidefinite and is entrywise non-negative.  DNN-based relaxations are a prominent family of semidefinite programming (SDP)-based relaxations used for solving non-convex quadratic optimization problems~\cite{bomzeDeKlerk2002,houTangToh2025rinnalplus}.  To address the second difficulty, we replace the inner maximization with its Lagrangian dual with the intention of deriving a single minimization instance.  However, as we later discuss, the resulting DNN relaxation does not satisfy Slater's condition, without which it is unclear if strong duality holds.  To address this, we apply a facial reduction step to obtain a reduced formulation whereby strong duality holds (whenever the inner problem is feasible).  These steps will allow us to obtain a single optimization instance in which we jointly minimize over the original decision $\theta \in \Theta$ together the dual variables associated with the inner maximization problem.  In particular, the formulation we obtain serves as an upper bound for \eqref{eq:special-cases-gwdro} whereby any gap that arises is simply due to the DNN relaxation.

The proofs of the results in this section are given in Appendix~\ref{app:proofs-dnn}.

In the remainder of this section, we fix a decision $\theta\in\Theta$ and we write $f$ in place of $f_{\theta}$.  We re-use the notation introduced in Section~\ref{subsec:problem-formulation}. %
In the following sequences of steps, we express the inner maximization problem within \eqref{eq:special-cases-gwdro} in matrix form, which will faciliate our derivation of the associated DNN relaxation.  In the remainder of this section, we let 
$$
z:=\operatorname{vec}(\pi)\in\mathbb{R}^{mn}
$$ 
denote the column-wise vectorization of $\pi$.  Suppose we index the rows and columns of $L\in\mathbb{R}^{mn\times mn}$ by pairs $(i,j)\in[m]\times[n]$.  Define the cost matrix
\begin{equation*}
    L_{(i,j),(k,l)}
    :=\bigl(D^{\mathcal{X}}_{ik}-D^{\mathcal{Y}}_{jl}\bigr)^2,
    \qquad i,k\in[m],\quad j,l\in[n]
\end{equation*}
so that the GW objective \eqref{eq:gw-distortion} can be written as 
$$
\mathcal{D}(\pi)=z^{\top}Lz.
$$
We proceed to the constraints.  First, we define the matrices
\begin{equation}
    A:=I_n\otimes\mathbf{1}_m^{\top}
    \in\mathbb{R}^{n\times mn},
    \qquad
    B:=\mathbf{1}_n^{\top}\otimes I_m
    \in\mathbb{R}^{m\times mn},
\label{eq:marginal-matrices}
\end{equation}
where $\otimes$ denotes the Kronecker product and $I_r$ is the
$r\times r$ identity matrix.  We then have $Az=\pi^{\top}\mathbf{1}_m = \widehat{\mu}_{\mathcal{Y}}$ and $Bz=\pi\mathbf{1}_n$.  Next, let $\boldsymbol{f}=(f(x_1),\ldots,f(x_m))^{\top}$ and define
\begin{equation*}
    c:=B^{\top}\boldsymbol{f}
    =\mathbf{1}_n\otimes\boldsymbol{f}.
\end{equation*}
Then the expected cost can be expressed as  $\boldsymbol{f}^{\top}Bz=c^{\top}z$.  Combining these, the inner maximization in \eqref{eq:special-cases-gwdro} can be equivalently expressed in the following form in which the decision variable is the vectorized matrix variable $z\in\mathbb{R}^{mn}$
\begin{equation}
\begin{aligned}
    \vartheta
    :=
    \max_{z\in\mathbb{R}^{mn}}
    \quad
    &c^{\top}z\\
    \text{s.t.}\quad
    &z^{\top}Lz\leq\rho^2,\\
    &Az=\widehat{\mu}_{\mathcal{Y}},\\
    &z\geq0.
\end{aligned}
\label{eq:worst-case-qcqp}
\end{equation}

\subsection{DNN relaxation}
\label{subsec:dnn-relaxation}

In this subsection, we consider solving the inner worst-case problem in \eqref{eq:special-cases-gwdro}.  As we noted earlier, the inner worst-case problem requires the solution of the GW problem, which is an instance of a non-convex QCQP, and is NP-Hard in general \cite{BGP:26,pardalosVavasis1991}.  To this end, a powerful family of methods for solving general QCQPs that involve binary constraints are those based on DNN-based relaxations~\cite{bomzeDeKlerk2002,houTangToh2025rinnalplus}.  These are powerful in the sense that they often return solutions that are close to being globally optimal in many practical settings.  In particular, within the specific context of the GW problem, prior experimental work in Chen et al.~\cite{chenNguyenKohSoh2024} as well as Tran et al. \cite{TNS:26} suggest that DNN-based relaxations are effective at recovering transportation plans that are globally optimal, or close to being globally optimal.  

As such, in the remainder of this section, we specialize the description of the DNN-based relaxation provided in~\cite{chenNguyenKohSoh2024} to our context of solving the inner worst-case problem (see also \cite{TNS:26}).  As we later note, a useful by-product of the DNN relaxation of the GW problem \cite{chenNguyenKohSoh2024} is that any solution $\hat{\pi}$ obtained from the relaxation is always a feasible coupling map.  In particular, by substituting $\hat{\pi}$ into the GW objective and checking if it coincides with the lower bound provided by the DNN relaxation, we have a simple test to check if $\hat{\pi}$ is also globally optimal.  We extend these ideas to our set-up where we optimize over a GW ball -- see Corollary~\ref{cor:dnn-exact} in the sequel.

\paragraph{Standard SDP relaxation.}
Suppose $Z \in \mathbb{S}^{mn}$ is a symmetric matrix and consider the following matrix
\begin{equation}
W
:=
\begin{bmatrix}
Z & z \\
z^{\top} & 1
\end{bmatrix}.
\label{eq:lifted-matrix}
\end{equation}
Then
\begin{equation*}
Z = zz^{\top}
\quad \Longleftrightarrow \quad
W \succeq 0
\quad \text{and} \quad
\operatorname{rank}(W) = 1.
\end{equation*}
In particular, \eqref{eq:worst-case-qcqp} is equivalent to
\begin{equation}
\vartheta
=
\sup_{\substack{z \in \mathbb{R}^{mn} \\ Z \in \mathbb{S}^{mn}}}
\left\{
c^{\top}z\;\middle|\;
\begin{aligned}
& \langle L,Z\rangle \leq \rho^2, \\
& Az = \widehat{\mu}_{\mathcal{Y}}, \quad z \geq 0, \\
& \begin{bmatrix} Z & z \\ z^{\top} & 1 \end{bmatrix} \succeq 0, \quad
\operatorname{rank}\!\begin{bmatrix} Z & z \\ z^{\top} & 1 \end{bmatrix} = 1
\end{aligned}
\right\}
\label{eq:rank-one-reformulation}
\end{equation}
The problem \eqref{eq:rank-one-reformulation} is non-convex due to the presence of the rank-one constraint.  By omitting the rank-one constraint we obtain the following relaxation, often referred to as the {\em Shor's relaxation}~\cite{shor1987quadratic}:
\begin{equation}
\vartheta_{\mathrm{SDP}}
:=
\sup_{\substack{z \in \mathbb{R}^{mn} \\ Z \in \mathbb{S}^{mn}}}
\left\{
c^{\top}z\;\middle|\;
\begin{aligned}
& \langle L,Z\rangle \leq \rho^2, \\
& Az = \widehat{\mu}_{\mathcal{Y}}, \quad z \geq 0, \\
& \begin{bmatrix} Z & z \\ z^{\top} & 1 \end{bmatrix} \succeq 0
\end{aligned}
\right\}
\label{eq:sdp-relaxation}
\end{equation}

\paragraph{Tightening the relaxation.}
The above SDP relaxation \eqref{eq:sdp-relaxation} can be further strengthened by adding valid linear conditions that are satisfied whenever $Z = zz^{\top}$.  First, by combining $Az = \widehat{\mu}_{\mathcal{Y}}$ with $Z = zz^{\top}$ we obtain
\begin{equation*}
AZ = Azz^{\top} = \widehat{\mu}_{\mathcal{Y}}z^{\top}.
\end{equation*}
Second, the conditions $z \geq 0$ and $Z = zz^{\top}$ imply $Z \geq 0$ componentwise.  In particular,  since $z$ and $Z$ are components of $W$ in \eqref{eq:lifted-matrix}, these non-negativity constraints imply $W \geq 0$.  Combining these, the matrix $W$ belongs to the doubly non-negative cone
\begin{equation*}
\mathcal{DNN}^{mn+1}
:=
\left\{
W \in \mathbb{S}^{mn+1} :
W \succeq 0,\ W \geq 0
\right\}.
\end{equation*}
By adding these valid constraints to \eqref{eq:sdp-relaxation} we obtain the following DNN relaxation
\begin{equation}
\vartheta_{\mathrm{DNN}}
:=
\sup_{\substack{z \in \mathbb{R}^{mn} \\ Z \in \mathbb{S}^{mn}}}
\left\{
c^{\top}z\;\middle|\;
\begin{aligned}
& \langle L,Z\rangle \leq \rho^2, \\
& Az = \widehat{\mu}_{\mathcal{Y}}, \quad AZ = \widehat{\mu}_{\mathcal{Y}}z^{\top}, \\
& \begin{bmatrix}
Z & z \\
z^{\top} & 1
\end{bmatrix} \in \mathcal{DNN}^{mn+1}
\end{aligned}
\right\}
\label{eq:dnn-relaxation}
\end{equation}

First, let's note the following ordering between the optimal values of the original problem \eqref{eq:worst-case-qcqp}, the standard SDP relaxation \eqref{eq:sdp-relaxation}, and the DNN relaxation \eqref{eq:dnn-relaxation}.

\begin{proposition}
\label{prop:dnn-sdp-bounds}
\begin{equation*}
\vartheta
\leq
\vartheta_{\mathrm{DNN}}
\leq
\vartheta_{\mathrm{SDP}}.
\end{equation*}
\end{proposition}

The DNN relaxation \eqref{eq:dnn-relaxation} provides a simple method to check if a given solution is the globally optimal solution to \eqref{eq:worst-case-qcqp}.

\begin{corollary}
\label{cor:dnn-exact}
Let $(\bar{z},\bar{Z})$ be an optimal solution to \eqref{eq:dnn-relaxation}. If $\bar{z}^{\top}L\bar{z} \leq \rho^2$, then
\begin{equation*}
\vartheta
=
\vartheta_{\mathrm{DNN}}
=
c^{\top}\bar{z}.
\end{equation*}
\end{corollary}

\subsection{Facial reduction}
\label{subsec:facial-reduction}

First, we show that Slater's condition fails for \eqref{eq:dnn-relaxation}. Slater's condition requires, in particular, the existence of a feasible solution that satisfies
\begin{equation*}
W=
\begin{bmatrix}
Z & z\\
z^{\top} & 1
\end{bmatrix}
\succ 0. 
\end{equation*}
Unfortunately, this condition cannot be satisfied because every $W$ that is feasible in \eqref{eq:dnn-relaxation} turns out to also be singular.

Recall that \eqref{eq:dnn-relaxation} imposes the constraints
\begin{equation}
Az = \widehat{\mu}_{\mathcal{Y}},
\qquad
AZ = \widehat{\mu}_{\mathcal{Y}}z^{\top}.
\label{eq:dnn-marginal-equalities}
\end{equation}
Define
\begin{equation*}
H
:=
\begin{bmatrix}
A & -\widehat{\mu}_{\mathcal{Y}}
\end{bmatrix}.
\end{equation*}
Direct multiplication gives
\begin{equation*}
HW
=
\begin{bmatrix}
AZ-\widehat{\mu}_{\mathcal{Y}}z^{\top}
&
Az-\widehat{\mu}_{\mathcal{Y}}
\end{bmatrix}.
\end{equation*}
Consequently, the equality constraints in \eqref{eq:dnn-marginal-equalities} are equivalent to
\begin{equation*}
HW = 0.
\end{equation*}
This implies $\operatorname{range}(W) \subseteq \ker(H)$.

The matrices $A = I_{n} \otimes \mathbf{1}_{m}^{\top}$ and $H$ have full row rank $n$. It follows that
\begin{equation*}
\dim\ker(H)
=
mn-n+1
<
mn+1.
\end{equation*}
Thus, $\ker(H)$ is a proper subspace of $\mathbb{R}^{mn+1}$, and every feasible $W$ satisfying $HW = 0$ satisfies
\begin{equation*}
\operatorname{rank}(W)
\leq
mn-n+1.
\end{equation*}
Therefore, every feasible $W$ is singular, so \eqref{eq:dnn-relaxation} does not satisfy Slater's condition.

Our next step is to apply a facial reduction step to obtain a reformulation that satisfies strong duality.  As a reminder, facial reduction is a method for reformulating a conic program in which one restricts the cone to a face of the cone containing all feasible solutions, without changing the feasible set or the optimal value~\cite{borweinWolkowicz1981}.  %
Our next step is to identify the corresponding face of the positive semidefinite cone. Since $W \succeq 0$, the equality $HW = 0$ is equivalent to
\begin{equation*}
\left\langle H^{\top}H,W\right\rangle
=
0.
\end{equation*}
Because $H^{\top}H$ is a nonzero positive semidefinite matrix, the set
\begin{equation*}
\left\{
W \in \mathbb{S}^{mn+1}_{+} :
\left\langle H^{\top}H,W\right\rangle = 0
\right\}
\end{equation*}
is a proper exposed face of $\mathbb{S}^{mn+1}_{+}$, where $\mathbb{S}^{r}_{+}$ denotes the cone of real symmetric positive semidefinite $r\times r$ matrices. Hence, the equality constraints in \eqref{eq:dnn-marginal-equalities} force every feasible matrix of \eqref{eq:dnn-relaxation} to lie within this proper face.  

Following the facial-reduction framework of Borwein and Wolkowicz~\cite{borweinWolkowicz1981}, we reformulate \eqref{eq:dnn-relaxation} over the lower-dimensional positive semidefinite cone corresponding to the face identified above.  More concretely, let $F \in \mathbb{R}^{mn \times (mn-n)}$ be a matrix whose columns form a basis of $\ker(A)$.  Choose $z^{0} \in \mathbb{R}^{mn}$ so that
\begin{equation}
Az^{0}
=
\widehat{\mu}_{\mathcal{Y}}.
\label{eq:facial-reduction-particular-point}
\end{equation}
When $m>1$, the solution to the linear system \eqref{eq:facial-reduction-particular-point} is not unique. Here, we choose
\begin{equation*}
z^{0}
=
\widehat{\mu}_{\mathcal{Y}}
\otimes
\frac{\mathbf{1}_{m}}{m},
\end{equation*}
so that $Az^{0}=\widehat{\mu}_{\mathcal{Y}}$.  Define
\begin{equation}
V
:=
\begin{bmatrix}
F & z^{0} \\
0^{\top} & 1
\end{bmatrix}.
\label{eq:facial-reduction-basis}
\end{equation}
The columns of $V$ form a basis of $\ker(H)$. It follows that
\begin{equation*}
W \succeq 0
\quad\text{and}\quad
HW = 0
\end{equation*}
if and only if there exists $Y \in \mathbb{S}^{mn-n+1}$ such that
\begin{equation*}
W
=
VYV^{\top},
\qquad
Y \succeq 0.
\end{equation*}
Next, we express the objective function and the GW constraint in terms of $Y$. Define
\begin{equation*}
\widetilde{C}
:=
\frac{1}{2}
\begin{bmatrix}
0 & c \\
c^{\top} & 0
\end{bmatrix},
\qquad
\widetilde{L}
:=
\begin{bmatrix}
L & 0 \\
0 & 0
\end{bmatrix}.
\end{equation*}
Then
\begin{equation*}
\left\langle \widetilde{C},W\right\rangle
=
c^{\top}z,
\qquad
\left\langle \widetilde{L},W\right\rangle
=
\left\langle L,Z\right\rangle.
\end{equation*}
Let
\begin{equation*}
\overline{C}:=V^{\top}\widetilde{C}V,
\qquad
\overline{L}:=V^{\top}\widetilde{L}V,
\qquad
e:=(0,\ldots,0,1)^{\top}
\in\mathbb{R}^{mn-n+1}.
\end{equation*}
By substituting $W = VYV^{\top}$ into \eqref{eq:dnn-relaxation} we obtain the equivalent optimization problem
\begin{equation}
\vartheta_{\mathrm{DNN}}
=
\sup_{Y \in \mathbb{S}^{mn-n+1}}
\left\{
\left\langle \overline{C},Y\right\rangle\;\middle|\;
\begin{aligned}
& \left\langle \overline{L},Y\right\rangle \leq \rho^2, \\
& \left\langle ee^{\top},Y\right\rangle = 1, \\
& VYV^{\top} \geq 0, \quad Y \succeq 0
\end{aligned}
\right\}
\label{eq:facially-reduced-dnn}
\end{equation}

The following theorem shows that the strong duality holds for the reduced problem \eqref{eq:facially-reduced-dnn}. This allows us to replace the inner maximization by a minimization over the dual variables. 
\begin{theorem}
\label{thm:facially-reduced-dnn-duality}
Suppose that \eqref{eq:facially-reduced-dnn} is feasible. The Lagrangian dual of \eqref{eq:facially-reduced-dnn} is given by
\begin{equation}
d_{\mathrm{DNN}}
:=
\inf_{\lambda,\tau,S}
\left\{
\lambda\rho^2+\tau\;\middle|\;
\begin{aligned}
& \lambda\overline{L}
  +\tau ee^{\top}
  -\overline{C}
  -V^{\top}SV
  \succeq0,\\
& \lambda\geq0,\qquad \tau\in\mathbb{R},\\
& S\in\mathbb{S}^{mn+1},\qquad S\geq0
\end{aligned}
\right\}
\label{eq:facially-reduced-dnn-dual}
\end{equation}
Moreover, strong duality holds in that
\begin{equation*}
\vartheta_{\mathrm{DNN}}=d_{\mathrm{DNN}}.
\end{equation*}
\end{theorem}

We now apply Theorem~\ref{thm:facially-reduced-dnn-duality} to replace the inner maximization of \eqref{eq:special-cases-gwdro} by its dual for each decision $\theta\in\Theta$. Let
$\boldsymbol{f}_{\theta}:=(f_{\theta}(x_1),\ldots,f_{\theta}(x_m))^{\top}\in\mathbb{R}^{m}$. Using $B$ from \eqref{eq:marginal-matrices} and $V$ from \eqref{eq:facial-reduction-basis}, define
\begin{equation*}
c_{\theta}:=B^{\top}\boldsymbol{f}_{\theta},
\qquad
\overline{C}_{\theta}
:=
\frac{1}{2}V^{\top}
\begin{bmatrix}
0&c_{\theta}\\
c_{\theta}^{\top}&0
\end{bmatrix}
V.
\end{equation*}
Denote by $\vartheta_{\mathrm{DNN}}(\theta)$ the optimal value of
\eqref{eq:facially-reduced-dnn} with this cost. Since the feasible set of
the worst-case problem does not depend on $\theta$, by applying
Theorem~\ref{thm:facially-reduced-dnn-duality} for each decision $\theta \in \Theta$ we obtain
\begin{equation}
\inf_{\theta\in\Theta}
\vartheta_{\mathrm{DNN}}(\theta)
=
\inf_{\theta,\lambda,\tau,S}
\left\{
\lambda\rho^2+\tau\;\middle|\;
\begin{aligned}
&\lambda\overline{L}
+\tau ee^{\top}
-\overline{C}_{\theta}
-V^{\top}SV
\succeq0,\\
&\theta\in\Theta,\qquad
\lambda\geq0,\qquad
\tau\in\mathbb{R},\\
&S\in\mathbb{S}^{mn+1},\qquad S\geq0
\end{aligned}
\right\}
\label{eq:gwdro-dnn-dual}
\end{equation}

The formulation \eqref{eq:gwdro-dnn-dual} replaces the DNN minimax problem
by a single minimization over $\theta$ and the dual variables. If
$\boldsymbol{f}_{\theta}$ is affine in $\theta$ and $\Theta$ has a
semidefinite representation, this is a semidefinite program that can be
solved using standard semidefinite programming solvers. Strong duality
preserves the value of the DNN relaxation. Since
$\vartheta(\theta)\leq\vartheta_{\mathrm{DNN}}(\theta)$ for every
$\theta$, the resulting optimal value is an upper bound on that of
\eqref{eq:gw-dro}, and equality with the original problem holds when the
DNN relaxation is tight for every decision.

\section{Numerical experiments}
\label{sec:numerical-experiments}

We investigate \eqref{eq:special-cases-gwdro} in two experiments: clustering over incomparable spaces, and generalized graph $k$-median. In both experiments, we replace the inner worst-case problem by its upper bound obtained from the dual of the doubly nonnegative relaxation given by \eqref{eq:facially-reduced-dnn-dual}, and the resulting optimization problem is \eqref{eq:gwdro-dnn-dual}.

\subsection{Clustering over incomparable spaces}
\label{subsec:spiral-clustering}

We consider the problem of clustering on the target space based on the distribution in the source space, where the target and source spaces are incomparable. Using the notation in Section~\ref{subsec:problem-formulation}, we are given the target metric space $(\mathcal{X},d_{\mathcal{X}})$ with $\mathcal{X}=\{x_i:i\in[m]\}$, the source metric space $(\mathcal{Y},d_{\mathcal{Y}})$ with $\mathcal{Y}=\{y_j:j\in[n]\}$, and an empirical source distribution $\widehat{\mu}_{\mathcal{Y},N}$ obtained from $N$ independent observations drawn from an underlying distribution $\mu_{\mathcal{Y}}$. To select $k$ centers, we consider solving the following instance of \eqref{eq:special-cases-gwdro}:
\begin{equation}
\label{eq:spiral-clustering-formulation}
\min_{\substack{C\subseteq\mathcal{X}\\|C|=k}}\;\max_{\mu\in\mathcal{U}_{\rho}^{\mathrm{GW}}(\widehat{\mu}_{\mathcal{Y},N})}\sum_{i\in[m]}\mu_i\min_{c\in C}d_{\mathcal{X}}(x_i,c)^2,
\end{equation}
where the ambiguity set is defined in \eqref{eq:finite-gw-ambiguity-set}. 

\begin{figure}[!htbp]
\centering
\includegraphics[width=\textwidth]{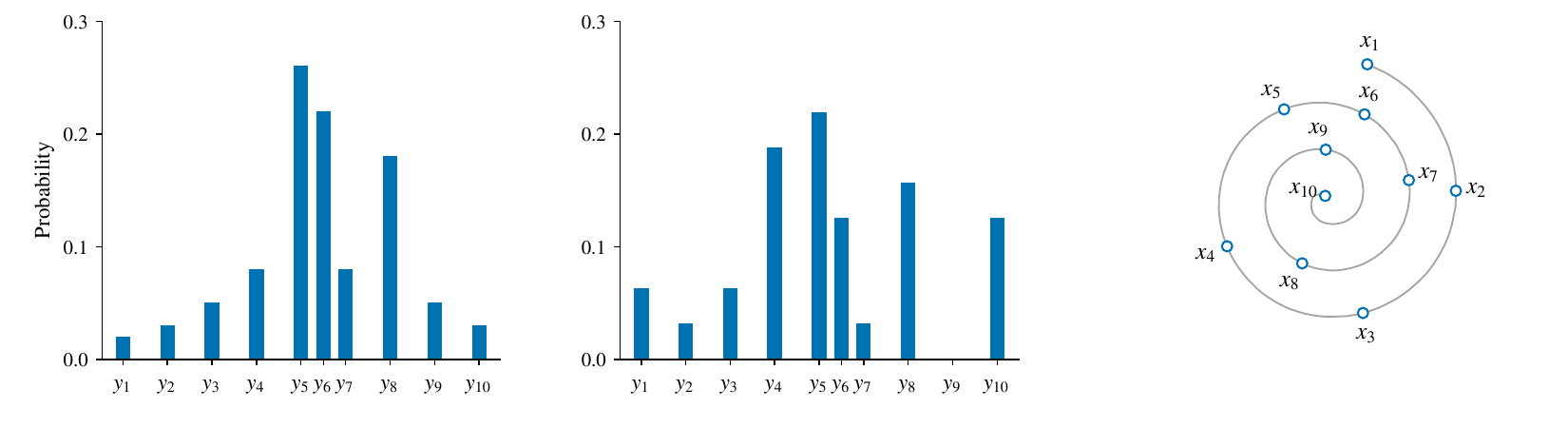}
\par\begin{subfigure}[t]{.32\textwidth}\caption{Population distribution.}\end{subfigure}\hfill
\begin{subfigure}[t]{.32\textwidth}\caption{Empirical distribution.}\end{subfigure}\hfill
\begin{subfigure}[t]{.32\textwidth}\caption{Target space.}\end{subfigure}
\caption{Source and target spaces for the clustering experiment. (a) Source population probabilities $\mu_{\mathcal{Y}}(\{y_i\})$. (b) Source empirical probabilities $\widehat{\mu}_{\mathcal{Y},32}(\{y_i\})$ in the first saved sample of $N=32$ observations. (c) Target points $x_1,\ldots,x_{10}$ on the spiral. The grey curve indicates the path used to measure arc length.}
\label{fig:spiral-data}
\end{figure}

In this experiment, the source space consists of $n=10$ points $y_1,\ldots,y_{10} \in \bbR$ on a line, and the target space consists of $m=10$ points $x_1,\ldots,x_{10} \in \bbR^2$ on a spiral, as shown in Figure~\ref{fig:spiral-data}. We select $k=3$ centers from the target points. The source metric $d_\mathcal{Y}$ measures distance along the line, and the target metric $d_\mathcal{X}$ measures arc length along the spiral. We index the points so that the normalized arc length of $x_i$ equals the position of $y_i$ on the line. For each sample size $N\in\{16,32,64,128,256,512\}$, we draw 30 independent samples from the source population distribution $\mu_{\mathcal{Y}}$ in Figure~\ref{fig:spiral-data}(a) to form the source empirical distribution $\widehat{\mu}_{\mathcal{Y},N}$.

In this experiment, the oracle is the center set $C^{\star}$ obtained by minimizing the expected loss under the target population distribution $\mu_{\mathcal{X}}$:
\begin{equation*}
C^{\star}\in\operatorname*{argmin}_{\substack{C\subseteq\mathcal{X}\\|C|=3}}\sum_{i\in[m]}\mu_{\mathcal{X},i}\min_{c\in C}d_{\mathcal{X}}(x_i,c)^2,
\end{equation*}
where the target population distribution $\mu_{\mathcal{X}}$ is defined by $\mu_{\mathcal{X}}(\{x_i\})=\mu_{\mathcal{Y}}(\{y_i\})$ for $i \in [10]$.

We compare GW-DRO in \eqref{eq:spiral-clustering-formulation} with srGW using \eqref{eq:special-cases-srgw-formulation}, with $\theta=C$, $\Theta=\{C\subseteq\mathcal{X}:|C|=3\}$, and $f_{\theta}(x_i)=f_C(x_i):=\min_{c\in C}d_{\mathcal{X}}(x_i,c)^2$. The estimate $\widehat{\mu}_{\mathcal{X}}^{\mathrm{srGW}}$ is obtained from \eqref{eq:special-cases-srgw-radius} with $\nu=\widehat{\mu}_{\mathcal{Y},N}$. For each source sample, we select the radius by five-fold cross-validation, following Mohajerin Esfahani and Kuhn~\cite[Section~7.2.B]{mohajerinesfahani2018data}. For both methods, we use Gurobi~\cite{gurobi} to solve the mixed-integer problems for center selection, and we use MOSEK~\cite{mosek} to solve the semidefinite subproblems given by the dual of the doubly nonnegative relaxation in \eqref{eq:facially-reduced-dnn-dual}.

Each method returns a center set $\widehat{C}$. We evaluate its performance by comparing $\widehat{C}$ with $C^{\star}$ using the following metrics:
\begin{enumerate}
\item \emph{Partition disagreement.} This measures the probability that a target point is assigned to a different cluster from the oracle, after matching the clusters to minimize disagreement:
\begin{equation*}
\min_{\substack{\sigma:\widehat{C}\to C^{\star}\\\sigma\text{ bijective}}}
\mathbb{P}_{x\sim\mu_{\mathcal{X}}}\!\left(
\sigma\!\left(\operatorname*{argmin}_{c\in\widehat{C}}d_{\mathcal{X}}(x,c)\right)
\neq\operatorname*{argmin}_{c\in C^{\star}}d_{\mathcal{X}}(x,c)
\right).
\end{equation*}
Each $\operatorname*{argmin}$ selects the center with the smallest index when several centers are equally near.
\item \emph{Excess population cost.} This quantifies the difference between the expected losses of $\widehat{C}$ and $C^{\star}$ under $\mu_{\mathcal{X}}$:
\begin{equation*}
\sum_{i\in[m]}\mu_{\mathcal{X},i}f_{\widehat{C}}(x_i)-\sum_{i\in[m]}\mu_{\mathcal{X},i}f_{C^{\star}}(x_i).
\end{equation*}
\end{enumerate}

\begin{figure}[!htbp]
\centering
\begin{subfigure}[t]{.48\textwidth}
\centering\includegraphics[width=\linewidth]{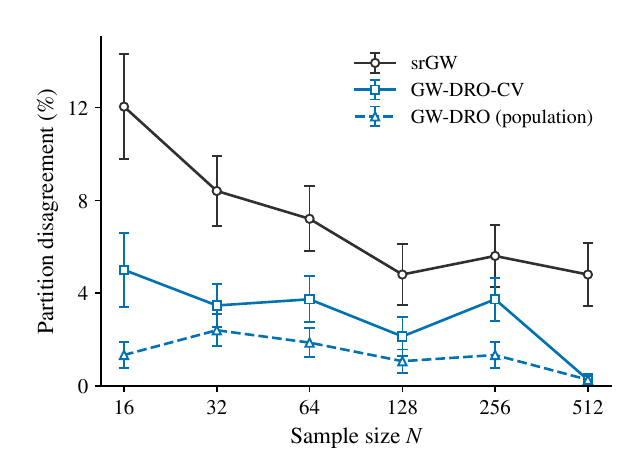}
\caption{Partition disagreement.}
\end{subfigure}\hfill
\begin{subfigure}[t]{.48\textwidth}
\centering\includegraphics[width=\linewidth]{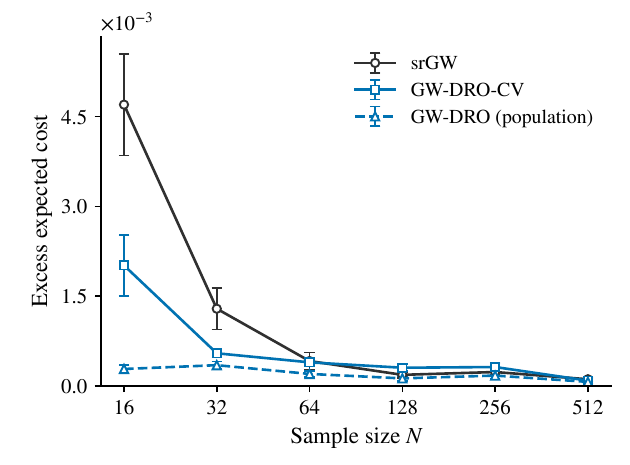}
\caption{Excess expected cost.}
\end{subfigure}
\caption{Comparison of GW-DRO and srGWon the clustering task over incomparable spaces. (a) Partition disagreement. (b) Excess expected cost. Both averaged over $30$ trails for each $N$. GW-DRO-CV denotes \eqref{eq:special-cases-gwdro} with the radius selected by cross-validation; GW-DRO (population) uses the population distribution to select the radius. }
\label{fig:spiral-results}
\end{figure}

The results are illustrated in Figure~\ref{fig:spiral-results}. We also include GW-DRO (population), which selects the radius in \eqref{eq:spiral-clustering-formulation} using the target population distribution, separately for each metric. Figure~\ref{fig:spiral-results}(a) shows an overall decrease in partition disagreement for all methods as the sample size increases. Averaged across the six sample sizes, GW-DRO-CV roughly halves the mean partition disagreement relative to srGW. Figure~\ref{fig:spiral-results}(b) shows a similar overall decrease in excess population cost, with the differences between the methods becoming smaller at larger sample sizes.

\FloatBarrier
\subsection{Stochastic graph $k$-median}
\label{subsec:branch-permutation-k-median}

We consider the stochastic $k$-median problem~\cite{anthony2010plant} on a connected undirected graph with vertex set $V=\{v_i:i\in[|V|]\}$ and edge set $E = \{\{v_i,v_j\}: v_i,v_j \in V\}$ with metric $d:V\times V \rightarrow \bbR_{+}$.  Let $p_v$ be the probability of demand at vertex $v$, with $p_v\geq0$ and $\sum_{v\in V}p_v=1$. For a nonempty center set $C\subseteq V$, recall from \eqref{eq:special-cases-outlier-center-distance} that $d(v,C):=\min_{c\in C}d(v,c)$ is the distance from $v$ to its nearest center. For an integer $1\leq k\leq |V|$, we solve
\begin{equation}
\label{eq:branch-population-cost}
\min_{\substack{C\subseteq V\\|C|=k}}\;\sum_{v\in V}p_v\,d(v,C).
\end{equation}

\begin{figure}[!htbp]
\centering
\begin{subfigure}[t]{.46\textwidth}
\centering\includegraphics[width=\linewidth]{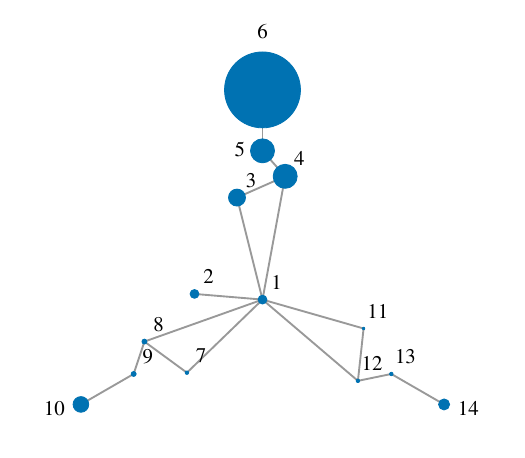}
\caption{Population distribution $p$.}
\end{subfigure}\hfill
\begin{subfigure}[t]{.46\textwidth}
\centering\includegraphics[width=\linewidth]{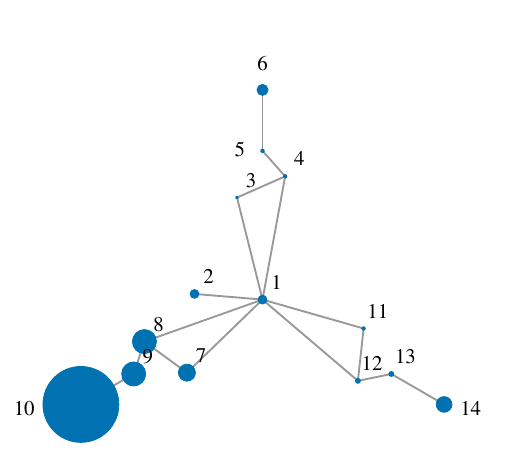}
\caption{Adversarial distribution $q$.}
\end{subfigure}
\caption{Visualization of graphs. The graph with $|V|=14$ vertices and vertex size proportional to probabilities. The adversarial distribution $q$ is obtained by redistributing the probability weights of $p$ among the branches. Edges indicate adjacency; their lengths in the drawing do not represent the metric $d$.}
\label{fig:branch-graph}
\end{figure}

In this experiment, we consider the graph in Figure~\ref{fig:branch-graph}, which has $|V|=14$ vertices and three branches with slightly different edge lengths. Our goal is to select $k=3$ centers on this graph. The population distribution $p$ places most of its mass on one branch.

First, we investigate the robustness of different methods to an adversarial perturbation of the observed distribution $p$. In the experiment in Section~\ref{subsec:spiral-clustering}, we consider uncertainty arising from sampling error in the empirical source distribution. Here, the observations are drawn from an adversarial distribution $q$ obtained by redistributing the probability weights of $p$ among the branches through vertex permutations. This perturbation approximately preserves isomorphism. We repeat the experiment independently ten times. In each repetition, we draw 500 vertices independently from $p$. For each sample size $N\in\{100,200,300,400,500\}$, we use the first $N$ observations and apply each of the six branch permutations, including the identity.  Each method receives only the observed vertex counts and the distance $d$.

We compare \eqref{eq:special-cases-gwdro} with empirical optimization and WDRO. For all three methods, $\theta=C$, $\Theta=\{C\subseteq V:|C|=3\}$, and $f_{\theta}(v)=f_C(v):=d(v,C)$. In \eqref{eq:special-cases-gwdro}, we take $\mathcal{X}=\mathcal{Y}=V$, $d_{\mathcal{X}}=d_{\mathcal{Y}}=d$, and source distribution $\widehat{\mu}_{\mathcal{Y},N}$. Empirical optimization solves \eqref{eq:branch-population-cost} with $p$ replaced by $\widehat{\mu}_{\mathcal{Y},N}$. WDRO solves \eqref{eq:dro} with the ambiguity set \eqref{eq:wdro-ambiguity-set} centered at $\widehat{\mu}_{\mathcal{Y},N}$ and transport cost $d(v,w)^2$ for $v,w\in V$. In this experiment, the oracle is a center set $C^{\star}$ obtained by solving \eqref{eq:branch-population-cost} with $k=3$ under the population distribution $p$. As in Section~\ref{subsec:spiral-clustering}, we use Gurobi~\cite{gurobi} to solve the mixed-integer problems for center selection in \eqref{eq:special-cases-gwdro}, and we use MOSEK~\cite{mosek} and SCS~\cite{scs} to solve the SDP subproblems in \eqref{eq:facially-reduced-dnn-dual}. For both \eqref{eq:special-cases-gwdro} and WDRO, we select the radius by three-fold cross-validation using the cost of the held-out observations at their observed vertices, and then fit the centers using all $N$ observations.  For comparison, we also include the performance of methods that select the radius to minimize the expected cost under the population $p$. 

Each method returns a center set $\widehat{C}$. We evaluate its cost under the population distribution $p$, that is, $\sum_{v\in V}p_v d(v,\widehat{C})$, and compare the oracle cost $\sum_{v\in V}p_v d(v,C^{\star})$, where $C^{\star}$ is a center set that minimizes the expected cost under $p$.

\begin{figure}[!htbp]
\centering
\begin{subfigure}[t]{.48\textwidth}
\centering\includegraphics[width=\linewidth]{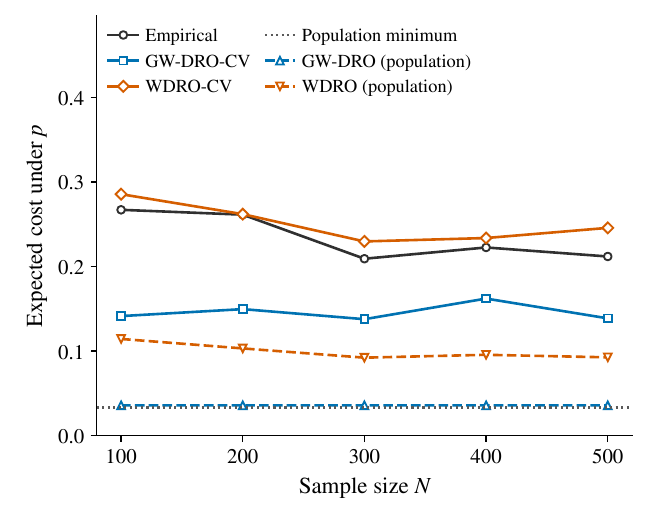}
\caption{Cost as the sample size varies.}
\end{subfigure}\hfill
\begin{subfigure}[t]{.48\textwidth}
\centering\includegraphics[width=\linewidth]{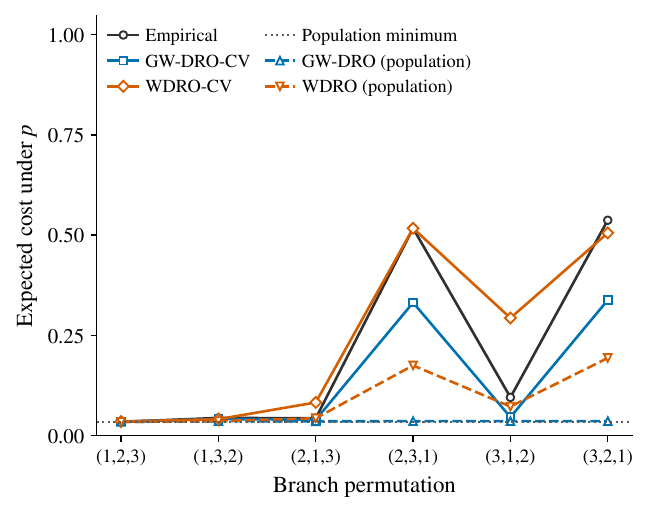}
\caption{Cost for each permutation at $N=500$.}
\end{subfigure}
\caption{Comparison among empirical optimization, \eqref{eq:special-cases-gwdro}, and WDRO on stochastic graph $k$-median. Population cost under $p$ is averaged over ten independent repetitions. Panel (a) averages over the six permutations at each sample size. Horizontal labels in panel (b) give the destinations of the three branches. Empirical denotes empirical optimization. GW-DRO-CV and WDRO-CV select the radius by cross-validation for \eqref{eq:special-cases-gwdro} and WDRO, respectively. The dashed curves labeled ``population'' use the population $p$ when selecting the radius. Dotted lines give the minimum population cost. }
\label{fig:branch-population-cost}
\end{figure}

The results are illustrated in Figure~\ref{fig:branch-population-cost}. Figure~\ref{fig:branch-population-cost}(a) shows that \eqref{eq:special-cases-gwdro} with cross-validation has lower mean population cost than empirical optimization and WDRO at every sample size considered. At $N=500$, the reduction is approximately $35\%$ and $44\%$, respectively. Figure~\ref{fig:branch-population-cost}(b) shows that its population cost also varies across permutations, with the largest reductions relative to empirical optimization occurring for permutations under which empirical optimization has the highest cost. Empirical optimization has slightly lower cost when no permutation is applied. We average the results over ten independent repetitions on the same graph. 

Second, we verify the relation between \eqref{eq:special-cases-gwdro} and robust graph $k$-median in \eqref{eq:special-cases-robust-k-median}. Suppose that $S_i = \{v_i\}$.  For every sufficiently large radius $\rho\geq\rho(\widehat{\mu}_{\mathcal{Y},N})$, where $\rho(\widehat{\mu}_{\mathcal{Y},N})$ is defined in \eqref{eq:special-cases-source-radius}, Proposition \ref{prop:special-cases-robust-k-median} gives
\begin{equation}
\label{eq:graph14-robust-equivalence}
\min_{\substack{C\subseteq V\\|C|=3}}\;
\max_{\mu\in\mathcal{U}_{\rho}^{\mathrm{GW}}}\sum_{v\in V}\mu_v d(v,C)
=\min_{\substack{C\subseteq V\\|C|=3}}\;\max_{v\in V}d(v,C).
\end{equation}
In this experiment, we verify numerically that \eqref{eq:graph14-robust-equivalence} holds, and show that the threshold of radius $\rho(\widehat{\mu}_{\mathcal{Y},N})$ for this to hold is tight. We solve the right-hand side of \eqref{eq:graph14-robust-equivalence} by evaluating all center sets and compare its optimal value with the computed objective of \eqref{eq:special-cases-gwdro} as the radius varies. We use the first sample with $N=500$ and the perturbation illustrated in Figure~\ref{fig:branch-graph}.

\begin{figure}[htbp!]
\centering
\includegraphics[width=0.8\textwidth]{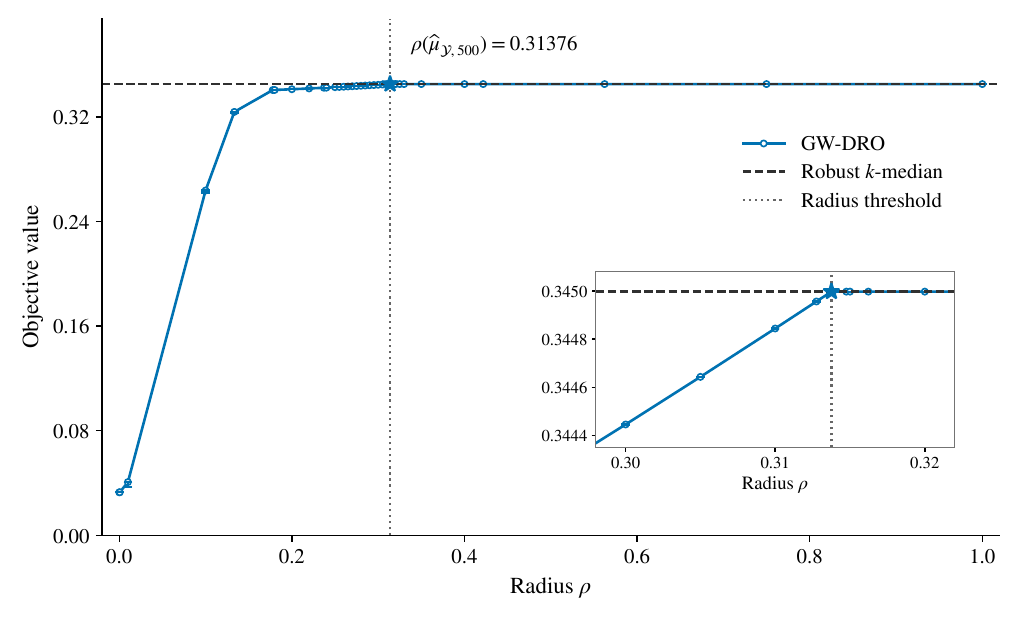}
\caption{The relation between stochastic graph $k$-median using \eqref{eq:special-cases-gwdro} and robust graph $k$-median. The optimal value of \eqref{eq:special-cases-gwdro} is nondecreasing in the radius $\rho$ and equals the optimal value of robust graph $k$-median when $\rho\geq\rho(\widehat{\mu}_{\mathcal{Y},500})$. The horizontal dashed line gives the robust graph $k$-median value, and the vertical dotted line denotes the radius threshold. Points give computed upper bounds, and vertical bars span the corresponding lower and upper bounds. The star marks the evaluation at the radius threshold, and the subfigure zooms in near the threshold.}
\label{fig:graph14-radius-objective}
\end{figure}

The results are illustrated in Figure~\ref{fig:graph14-radius-objective}. The computed objective of \eqref{eq:special-cases-gwdro} increases with the radius until it reaches the optimal value of robust graph $k$-median. At every evaluated radius satisfying $\rho\geq\rho(\widehat{\mu}_{\mathcal{Y},500})$, the two values agree, consistent with Proposition~\ref{prop:special-cases-robust-k-median}. The subfigure indicates that $\rho(\widehat{\mu}_{\mathcal{Y},500})$ is close to the smallest radius at which the two optimal values coincide for this dataset, indicating this threshold is tight.

\FloatBarrier

\bibliography{biblio}

\clearpage
\appendix
\section{Proofs}
\label{app:proofs}

\subsection{Proofs for Section~\ref{sec:gwdro-special-cases}}
\label{app:proofs-examples}

\begin{proof}[Proof of Proposition~\ref{prop:special-cases-srgw}]
It suffices to show that the GW ambiguity set at radius $\rho_{\min}(\widehat{\mu}_{\mathcal{Y}})$ contains only $\widehat{\mu}_{\mathcal{X}}^{\mathrm{srGW}}$.
By the definition of $\rho_{\min}$ in \eqref{eq:special-cases-srgw-radius},
\begin{equation*}
\operatorname{GW}_2^2\bigl(\mu,\widehat{\mu}_{\mathcal{Y}}\bigr)\geq\rho_{\min}(\widehat{\mu}_{\mathcal{Y}})^2
\qquad(\mu\in\Delta_m).
\end{equation*}
Combining this inequality with the definition of $\mathcal{U}_{\rho_{\min}(\widehat{\mu}_{\mathcal{Y}})}^{\mathrm{GW}}(\widehat{\mu}_{\mathcal{Y}})$ in \eqref{eq:finite-gw-ambiguity-set} and the uniqueness of the minimizer gives
\begin{equation*}
\mathcal{U}_{\rho_{\min}(\widehat{\mu}_{\mathcal{Y}})}^{\mathrm{GW}}(\widehat{\mu}_{\mathcal{Y}})
=\left\{\mu\in\Delta_m:
\operatorname{GW}_2^2\bigl(\mu,\widehat{\mu}_{\mathcal{Y}}\bigr)=\rho_{\min}(\widehat{\mu}_{\mathcal{Y}})^2\right\}
=\left\{\widehat{\mu}_{\mathcal{X}}^{\mathrm{srGW}}\right\}.
\end{equation*}
It follows that, for every $\theta\in\Theta$,
\begin{equation*}
\max_{\mu\in\mathcal{U}_{\rho_{\min}(\widehat{\mu}_{\mathcal{Y}})}^{\mathrm{GW}}(\widehat{\mu}_{\mathcal{Y}})}\sum_{i\in[m]}\mu_i f_{\theta}(x_i)
=\sum_{i\in[m]}\widehat{\mu}_{\mathcal{X},i}^{\mathrm{srGW}}f_{\theta}(x_i).
\end{equation*}
Minimizing over $\theta\in\Theta$ proves that \eqref{eq:special-cases-gwdro} at radius $\rho_{\min}(\widehat{\mu}_{\mathcal{Y}})$ is equivalent to \eqref{eq:special-cases-srgw-formulation}.
\end{proof}

\begin{proof}[Proof of Proposition~\ref{prop:special-cases-cvar}]
By the definition of the GW ambiguity set in \eqref{eq:finite-gw-ambiguity-set}, using the source distribution and distance matrices in Proposition~\ref{prop:special-cases-cvar} and taking $\rho=0$, the GW ambiguity set is
\begin{equation}
\label{eq:special-cases-cvar-finite-ambiguity-set}
\mathcal{U}_{0}^{\mathrm{GW}}(\widehat{\mu}_{\mathcal{Y}})
=\left\{\pi\mathbf{1}_n:\pi\in\mathbb{R}_{+}^{m\times n},\,\pi^{\top}\mathbf{1}_m=\mathbf{1}_n/n,\,\mathcal{D}(\pi)=\sum_{\substack{i,j\in[m]\\k,l\in[n]}}\bigl(D^{\mathcal{X}}_{ij}-D^{\mathcal{Y}}_{kl}\bigr)^2\pi_{ik}\pi_{jl}\leq 0 \right\}.
\end{equation}
We start by showing that, under the assumptions of Proposition~\ref{prop:special-cases-cvar}, the GW ambiguity set is given by
\begin{equation}
\mathcal{U}_{0}^{\mathrm{GW}}(\widehat{\mu}_{\mathcal{Y}})
=\{\operatorname{Unif}(S):S\subseteq\mathcal{X},\ |S|=n\}.
\label{eq:special-cases-cvar-ambiguity-set}
\end{equation}
\noindent$(\subseteq)$. Let $\mu\in\mathcal{U}_{0}^{\mathrm{GW}}(\widehat{\mu}_{\mathcal{Y}})$. By the expression for the GW ambiguity set in \eqref{eq:special-cases-cvar-finite-ambiguity-set}, there exists a coupling $\pi\in\mathbb{R}_{+}^{m\times n}$ satisfying
\begin{equation*}
\sum_{i\in[m]}\pi_{ik}=\frac1n\quad(k\in[n]),
\qquad \sum_{k\in[n]}\pi_{ik}=\mu_i\quad(i\in[m]),
\qquad \mathcal{D}(\pi)\leq0.
\end{equation*}
For $D^\mathcal{X}$ and $D^\mathcal{Y}$ defined in the statement of the proposition, $(D^{\mathcal{X}}_{ij}-D^{\mathcal{Y}}_{kl})^2$ equals $1$ when exactly one of the equalities $i=j$ and $k=l$ holds, and equals $0$ otherwise. Substituting this into \eqref{eq:gw-distortion} gives
\begin{equation*}
\mathcal{D}(\pi)
=\sum_{k\in[n]}\sum_{\substack{i,j\in[m]\\i\neq j}}\pi_{ik}\pi_{jk}
+\sum_{i\in[m]}\sum_{\substack{k,l\in[n]\\k\neq l}}\pi_{ik}\pi_{il} \ge 0, 
\end{equation*}
where the last inequality follows since $\pi\geq0$. Since $\mathcal{D}(\pi)\leq0$, every summand must equal zero. Hence,
\begin{equation*}
\pi_{ik}\pi_{jk}=0\quad(i\neq j,\ k\in[n]),
\qquad
\pi_{ik}\pi_{il}=0\quad(k\neq l,\ i\in[m]).
\end{equation*}
The first equality implies that each column of $\pi$ contains at most one positive entry. Together with $\sum_{i\in[m]}\pi_{ik}=1/n$, this shows that column $k$ contains exactly one positive entry, equal to $1/n$. Let $i_k$ denote its row index of the nonzero entry in the $k$-th column of $\pi$. The second equality implies that each row contains at most one positive entry. In particular, if $i_k=i_l$ for some $k\neq l$, then
\begin{equation*}
\pi_{i_k k}\pi_{i_k l}=\frac{1}{n^2}>0,
\end{equation*}
a contradiction. Thus $i_k\neq i_l$ whenever $k\neq l$, and $S:=\{x_{i_k}:k\in[n]\}$ contains exactly $n$ distinct points. Consequently, the row sums of $\pi$ satisfy
\begin{equation*}
\mu_i=\sum_{k\in[n]}\pi_{ik}
=\begin{cases}
1/n,&x_i\in S,\\
0,&x_i\notin S.
\end{cases}
\end{equation*}
Consequently, $\mu=\operatorname{Unif}(S)$ with $|S|=n$, which proves the inclusion $\subseteq$.

\noindent$(\supseteq)$. Fix $S=\{x_{i_1},\ldots,x_{i_n}\}\subseteq\mathcal{X}$ with $n$ distinct points, and define
\begin{equation*}
\pi_{ik}=\begin{cases}
1/n,&i=i_k,\\
0,&i\neq i_k.
\end{cases}
\end{equation*}
Each column of $\pi$ sums to $1/n$. Since the indices $i_1,\ldots,i_n$ are distinct, the $i_k$-th row of $\pi$ sums to $1/n$ for every $k\in[n]$, and every other row of $\pi$ sums to zero. Hence $\pi$ is a coupling between $\operatorname{Unif}(S)$ and $\widehat{\mu}_{\mathcal{Y}}$. Moreover, substituting the constructed $\pi$ into \eqref{eq:gw-distortion} gives
\begin{equation*}
\mathcal{D}(\pi)
=\sum_{k,l\in[n]}\sum_{i,j\in[m]}\bigl(D^{\mathcal{X}}_{ij}-D^{\mathcal{Y}}_{kl}\bigr)^2\pi_{ik}\pi_{jl}\!
=\frac{1}{n^2}\sum_{k,l\in[n]}\sum_{\substack{i,j\in[m]\\i=i_k,\ j=i_l}}\bigl(D^{\mathcal{X}}_{ij}-D^{\mathcal{Y}}_{kl}\bigr)^2\!
=\frac{1}{n^2}\sum_{k,l\in[n]}\bigl(D^{\mathcal{X}}_{i_k i_l}-D^{\mathcal{Y}}_{kl}\bigr)^2
=0,
\end{equation*}
where the second equality uses $\pi_{ik}\pi_{jl}=1/n^2$ when $i=i_k$ and $j=i_l$, and $\pi_{ik}\pi_{jl}=0$ otherwise. The third equality evaluates the sum over $i,j$, and the last equality follows from $D^{\mathcal{X}}_{i_k i_l}=D^{\mathcal{Y}}_{kl}$ for all $k,l\in[n]$. Thus $\operatorname{Unif}(S)\in\mathcal{U}_{0}^{\mathrm{GW}}(\widehat{\mu}_{\mathcal{Y}})$, proving the inclusion $\supseteq$.

Combining Lemma~\ref{lem:special-cases-capped-simplex-hull} and \eqref{eq:special-cases-cvar-ambiguity-set} gives
\begin{equation}
\operatorname{conv}\bigl(\mathcal{U}_{0}^{\mathrm{GW}}(\widehat{\mu}_{\mathcal{Y}})\bigr)=\conv\bigl\{\operatorname{Unif}(S):\ S\subseteq\mathcal{X},\ |S|=n\bigr\}=\mathcal{C}_n.
\label{eq:special-cases-cvar-ambiguity-convex-hull}
\end{equation}
It follows that, for every $\theta\in\Theta$,
\begin{equation*}
\begin{aligned}
\max_{\mu\in\mathcal{U}_{0}^{\mathrm{GW}}(\widehat{\mu}_{\mathcal{Y}})}\sum_{i\in[m]}\mu_i f_{\theta}(x_i)
&=\max_{\mu\in\operatorname{conv}(\mathcal{U}_{0}^{\mathrm{GW}}(\widehat{\mu}_{\mathcal{Y}}))}\sum_{i\in[m]}\mu_i f_{\theta}(x_i)\\
&=\max_{\mu\in\mathcal{C}_n}\sum_{i\in[m]}\mu_i f_{\theta}(x_i)
=\operatorname{CVaR}_{1-n/m}\bigl(f_{\theta}(x)\bigr),
\end{aligned}
\end{equation*}
where the first equality holds because the objective is linear in $\mu$, the second equality follows from \eqref{eq:special-cases-cvar-ambiguity-convex-hull}, and the third equality follows from Lemma~\ref{lem:special-cases-cvar-dual}. Minimizing over $\theta\in\Theta$ proves that \eqref{eq:special-cases-gwdro} is equivalent to \eqref{eq:special-cases-cvar-formulation} with $\alpha=1-n/m$.
\end{proof}

\begin{proof}[Proof of Proposition~\ref{prop:special-cases-group-dro}]
The proof is similar to the proof of Proposition~\ref{prop:special-cases-cvar}. Recall from the definition of the GW ambiguity set in \eqref{eq:finite-gw-ambiguity-set} that, under the assumptions of Proposition~\ref{prop:special-cases-group-dro},
\begin{equation}
\label{eq:special-cases-group-dro-finite-ambiguity-set}
\mathcal{U}_{0}^{\mathrm{GW}}(\widehat{\mu}_{\mathcal{Y}})
=\left\{\pi\mathbf{1}_n:\pi\in\mathbb{R}_{+}^{m\times n},\,\pi^{\top}\mathbf{1}_m=\mathbf{1}_n/n,\,\mathcal{D}(\pi)=\sum_{\substack{i,j\in[m]\\k,l\in[n]}}\bigl(D^{\mathcal{X}}_{ij}-D^{\mathcal{Y}}_{kl}\bigr)^2\pi_{ik}\pi_{jl}\leq 0\right\}.
\end{equation}
We start by proving the following equality:
\begin{equation}
\mathcal{U}_{0}^{\mathrm{GW}}(\widehat{\mu}_{\mathcal{Y}})
=\{\operatorname{Unif}(\mathcal{X}_g):g\in[G]\}.
\label{eq:special-cases-group-dro-ambiguity-set}
\end{equation}

\noindent$(\subseteq)$. Let $\mu\in\mathcal{U}_{0}^{\mathrm{GW}}(\widehat{\mu}_{\mathcal{Y}})$. By the expression for the GW ambiguity set in \eqref{eq:special-cases-group-dro-finite-ambiguity-set}, there exists a coupling $\pi\in\mathbb{R}_{+}^{m\times n}$ satisfying
\begin{equation*}
\sum_{i\in[m]}\pi_{ik}=\frac1n\quad(k\in[n]),
\qquad \sum_{k\in[n]}\pi_{ik}=\mu_i\quad(i\in[m]),
\qquad \mathcal{D}(\pi)\leq0.
\end{equation*}
Substituting the target and source distances into \eqref{eq:gw-distortion} gives
\begin{equation*}
\mathcal{D}(\pi)
=\sum_{k\in[n]}\sum_{\substack{i,j\in[m]\\i\neq j}}\bigl(D^{\mathcal{X}}_{ij}\bigr)^2\pi_{ik}\pi_{jk}
+\sum_{i\in[m]}\sum_{\substack{k,l\in[n]\\k\neq l}}\pi_{ik}\pi_{il}
+\sum_{\substack{k,l\in[n]\\k\neq l}}\ \sum_{\substack{i,j\in[m]\\x_i,x_j\text{ in different groups}}}\pi_{ik}\pi_{jl} \geq 0,
\end{equation*}
where the last inequality follows since $\pi\geq0$. Since $\mathcal{D}(\pi)\leq0$, every summand must equal zero. Moreover, $D^{\mathcal{X}}_{ij}\geq1$ for $i\neq j$. Hence, for $i,j\in[m]$ and $k,l\in[n]$,
\begin{equation*}
\pi_{ik}\pi_{jk}=0\quad(i\neq j),
\qquad \pi_{ik}\pi_{il}=0\quad(k\neq l),
\qquad \pi_{ik}\pi_{jl}=0\quad(k\neq l,\ x_i,x_j\text{ in different groups}).
\end{equation*}
The first equality implies that each column of $\pi$ contains at most one positive entry. Together with $\sum_{i\in[m]}\pi_{ik}=1/n$, this shows that column $k$ of $\pi$ contains exactly one positive entry, equal to $1/n$. Let $i_k$ denote the row index of the nonzero entry in the $k$-th column of $\pi$. The second equality implies that each row of $\pi$ contains at most one positive entry. In particular, if $i_k=i_l$ for some $k\neq l$, then
\begin{equation*}
\pi_{i_k k}\pi_{i_k l}=\frac{1}{n^2}>0,
\end{equation*}
a contradiction. Thus $i_k\neq i_l$ whenever $k\neq l$, and the $n$ points $x_{i_1},\ldots,x_{i_n}$ are distinct. The third equality implies that $x_{i_k}$ and $x_{i_l}$ lie in the same group for $k\neq l$, since otherwise $\pi_{i_k k}\pi_{i_l l}=1/n^2>0$. Thus $x_{i_1},\ldots,x_{i_n}$ lie in a single group $\mathcal{X}_g$. Since $|\mathcal{X}_g|=n$, these points comprise all of $\mathcal{X}_g$. Consequently, the row sums of $\pi$ satisfy
\begin{equation*}
\mu_i=\sum_{k\in[n]}\pi_{ik}
=\begin{cases}
1/n,&x_i\in\mathcal{X}_g,\\
0,&x_i\notin\mathcal{X}_g.
\end{cases}
\end{equation*}
Consequently, $\mu=\operatorname{Unif}(\mathcal{X}_g)$ for some $g\in[G]$, which proves the inclusion $\subseteq$.

\noindent$(\supseteq)$. Fix $g\in[G]$ and write $\mathcal{X}_g=\{x_{i_1},\ldots,x_{i_n}\}$. We prove that $\operatorname{Unif}(\mathcal{X}_g)\in\mathcal{U}_{0}^{\mathrm{GW}}(\widehat{\mu}_{\mathcal{Y}})$. To this end, define
\begin{equation*}
\pi_{ik}=\begin{cases}
1/n,&i=i_k,\\
0,&i\neq i_k.
\end{cases}
\end{equation*}
Each column of $\pi$ sums to $1/n$. Since the indices $i_1,\ldots,i_n$ are distinct, the $i_k$-th row of $\pi$ sums to $1/n$ for every $k\in[n]$, and every other row of $\pi$ sums to zero. Hence $\pi$ is a coupling between $\operatorname{Unif}(\mathcal{X}_g)$ and $\widehat{\mu}_{\mathcal{Y}}$. Moreover, substituting the constructed $\pi$ into \eqref{eq:gw-distortion} gives
\begin{equation*}
\mathcal{D}(\pi)
=\sum_{k,l\in[n]}\sum_{i,j\in[m]}\bigl(D^{\mathcal{X}}_{ij}-D^{\mathcal{Y}}_{kl}\bigr)^2\pi_{ik}\pi_{jl}\!
=\frac{1}{n^2}\sum_{k,l\in[n]}\sum_{\substack{i,j\in[m]\\i=i_k,\ j=i_l}}\bigl(D^{\mathcal{X}}_{ij}-D^{\mathcal{Y}}_{kl}\bigr)^2\!
=\frac{1}{n^2}\sum_{k,l\in[n]}\bigl(D^{\mathcal{X}}_{i_k i_l}-D^{\mathcal{Y}}_{kl}\bigr)^2
=0,
\end{equation*}
where the second equality uses $\pi_{ik}\pi_{jl}=1/n^2$ when $i=i_k$ and $j=i_l$, and $\pi_{ik}\pi_{jl}=0$ otherwise. The third equality evaluates the sum over $i,j$, and the last equality follows from $D^{\mathcal{X}}_{i_k i_l}=D^{\mathcal{Y}}_{kl}$ for all $k,l\in[n]$, since $x_{i_1},\ldots,x_{i_n}$ are distinct points in the same group. Thus $\operatorname{Unif}(\mathcal{X}_g)\in\mathcal{U}_{0}^{\mathrm{GW}}(\widehat{\mu}_{\mathcal{Y}})$, proving the inclusion $\supseteq$.

It follows from \eqref{eq:special-cases-group-dro-ambiguity-set} that, for every $\theta\in\Theta$,
\begin{equation*}
\max_{\mu\in\mathcal{U}_{0}^{\mathrm{GW}}(\widehat{\mu}_{\mathcal{Y}})}\sum_{i\in[m]}\mu_i f_{\theta}(x_i)
=\max_{g\in[G]}\frac1n\sum_{x\in\mathcal{X}_g}f_{\theta}(x).
\end{equation*}
Minimizing over $\theta\in\Theta$ proves that \eqref{eq:special-cases-gwdro} is equivalent to \eqref{eq:special-cases-group-dro-formulation}.
\end{proof}

\begin{proof}[Proof of Proposition~\ref{prop:special-cases-robust-k-median}]
By the definition of the GW ambiguity set in \eqref{eq:finite-gw-ambiguity-set}, the GW ambiguity set for the metric spaces in Proposition~\ref{prop:special-cases-robust-k-median} is
\begin{equation*}
\mathcal{U}_{\rho}^{\mathrm{GW}}(\widehat{\mu}_{\mathcal{Y}})
=\left\{\pi\mathbf{1}_n:\pi\in\mathbb{R}_{+}^{G\times n},\,\pi^{\top}\mathbf{1}_G=\widehat{\mu}_{\mathcal{Y}},\,\sum_{i,j\in[G]}\sum_{k,l\in[n]}\bigl(d_{\mathcal{X}}(S_i,S_j)-d_{\mathcal{Y}}(y_k,y_l)\bigr)^2\pi_{ik}\pi_{jl}\leq\rho^2\right\}.
\end{equation*}
We first show that the GW ambiguity set in \eqref{eq:special-cases-gwdro} contains every distribution concentrated on one of the sets $S_1,\ldots,S_G$. That is,
\begin{equation*}
\{e_i:i\in[G]\}\subseteq\mathcal{U}_{\rho}^{\mathrm{GW}}(\widehat{\mu}_{\mathcal{Y}}),
\end{equation*}
where $e_i\in\mathbb{R}^{G}$ has entry $1$ at index $i$ and entry $0$ elsewhere, and represents the distribution assigning probability $1$ to $S_i$.
Fix $i\in[G]$ and consider the following construction of $\pi\in\mathbb{R}_{+}^{G\times n}$: 
\begin{equation}
\label{eq:special-cases-point-mass-coupling}
\pi_{jl}=\begin{cases}
\widehat{\mu}_{\mathcal{Y},l},&j=i,\\
0,&j\neq i,
\end{cases}
\qquad j\in[G],\quad l\in[n].
\end{equation}
Since $\widehat{\mu}_{\mathcal{Y}}\in\Delta_n$, the source and target probabilities of this coupling satisfy $\pi^{\top}\mathbf{1}_{G}=\widehat{\mu}_{\mathcal{Y}}$ and $\pi\mathbf{1}_{n}=e_i$.  Moreover, substituting the constructed $\pi$ into $\mathcal{D}(\pi)$ defined in \eqref{eq:gw-distortion} gives
\begin{equation*}
\mathcal{D}(\pi)
=\sum_{j,j'\in[G]}\sum_{k,l\in[n]}\bigl(D^{\mathcal{X}}_{jj'}-D^{\mathcal{Y}}_{kl}\bigr)^2\pi_{jk}\pi_{j'l}
=\sum_{k,l\in[n]}\widehat{\mu}_{\mathcal{Y},k}\widehat{\mu}_{\mathcal{Y},l}\bigl(D^{\mathcal{X}}_{ii}-D^{\mathcal{Y}}_{kl}\bigr)^2
=\rho(\widehat{\mu}_{\mathcal{Y}})^2\leq\rho^2.
\end{equation*}
Here, the second equality follows from \eqref{eq:special-cases-point-mass-coupling}, since $\pi_{ik}=\widehat{\mu}_{\mathcal{Y},k}$ for $k\in[n]$ and $\pi_{jk}=0$ for $j\in[G]\setminus\{i\}$ and $k\in[n]$. The third equality follows from $D^{\mathcal{X}}_{ii}=d_{\mathcal{X}}(S_i,S_i)=0$ and the definition of $\rho(\widehat{\mu}_{\mathcal{Y}})$ in \eqref{eq:special-cases-source-radius}. The inequality follows from the assumption on $\rho$ in Proposition~\ref{prop:special-cases-robust-k-median}. Thus $e_i\in\mathcal{U}_{\rho}^{\mathrm{GW}}(\widehat{\mu}_{\mathcal{Y}})$ for every $i\in[G]$.

We next show that, for every $\theta\in\Theta$,
\begin{equation}
\label{eq:special-cases-robust-inner-equality}
\max_{\mu\in\mathcal{U}_{\rho}^{\mathrm{GW}}(\widehat{\mu}_{\mathcal{Y}})}\sum_{i\in[G]}\mu_i f_{\theta}(S_i)
=\max_{i\in[G]}f_{\theta}(S_i).
\end{equation}
Fix $\theta\in\Theta$. We first show that the left-hand side of \eqref{eq:special-cases-robust-inner-equality} is at least the right-hand side. Since $e_i\in\mathcal{U}_{\rho}^{\mathrm{GW}}(\widehat{\mu}_{\mathcal{Y}})$ for every $i\in[G]$,
\begin{equation*}
\max_{i\in[G]}f_{\theta}(S_i)
=\max_{i\in[G]}\sum_{j\in[G]}(e_i)_j f_{\theta}(S_j)
\leq\max_{\mu\in\mathcal{U}_{\rho}^{\mathrm{GW}}(\widehat{\mu}_{\mathcal{Y}})}\sum_{j\in[G]}\mu_j f_{\theta}(S_j).
\end{equation*}
We next show that the left-hand side of \eqref{eq:special-cases-robust-inner-equality} is at most the right-hand side. Every $\mu\in\mathcal{U}_{\rho}^{\mathrm{GW}}(\widehat{\mu}_{\mathcal{Y}})$ satisfies $\mu_j\geq0$ for $j\in[G]$ and $\sum_{j\in[G]}\mu_j=1$. Therefore,
\begin{equation*}
\sum_{j\in[G]}\mu_j f_{\theta}(S_j)
\leq\sum_{j\in[G]}\mu_j\max_{i\in[G]}f_{\theta}(S_i)
=\max_{i\in[G]}f_{\theta}(S_i).
\end{equation*}
Taking the maximum over $\mu\in\mathcal{U}_{\rho}^{\mathrm{GW}}(\widehat{\mu}_{\mathcal{Y}})$ establishes the inequality $\leq$ in \eqref{eq:special-cases-robust-inner-equality}. Minimizing \eqref{eq:special-cases-robust-inner-equality} over $\theta\in\Theta$ and using the definition of $f_{\theta}$ in Proposition~\ref{prop:special-cases-robust-k-median} proves that \eqref{eq:special-cases-gwdro} is equivalent to \eqref{eq:special-cases-robust-k-median}.
\end{proof}

\begin{proof}[Proof of Lemma~\ref{lem:special-cases-outlier-duality}]
Fix a nonempty $C\subseteq V$. We start by reformulating the minimization on the left-hand side of \eqref{eq:special-cases-outlier-duality} as a linear program. We then show that the maximization on the right-hand side of \eqref{eq:special-cases-outlier-duality} is the dual of this linear program. 

First, We claim that the following equalities hold:
\begin{equation}
\label{eq:special-cases-outlier-retained-lp}
\min_{\substack{O\subseteq V\\|O|=z}}\sum_{v\in V\setminus O}d(v,C)
=\min_{\substack{w_v\in\{0,1\}\ (v\in V)\\\sum_{v\in V}w_v=|V|-z}}\sum_{v\in V}w_vd(v,C)
=\min_{\substack{w_v\in[0,1]\ (v\in V)\\\sum_{v\in V}w_v=|V|-z}}\sum_{v\in V}w_vd(v,C).
\end{equation}
To prove the first equality, let $O\subseteq V$ with $|O|=z$ and define $w=(w_v)_{v\in V}$ by
\begin{equation*}
w_v=\begin{cases}
0,&v\in O,\\
1,&v\in V\setminus O.
\end{cases}
\end{equation*}
The vector $w$ satisfies
\begin{equation*}
\sum_{v\in V}w_v=|V\setminus O|=|V|-z,\qquad
\sum_{v\in V}w_vd(v,C)=\sum_{v\in V\setminus O}d(v,C).
\end{equation*}
Thus, $w$ is feasible for the binary linear program in \eqref{eq:special-cases-outlier-retained-lp} with objective value $\sum_{v\in V\setminus O}d(v,C)$.
For any $w$ feasible for the binary linear program in \eqref{eq:special-cases-outlier-retained-lp}, let $O=\{v\in V:w_v=0\}$. Since $w_v\in\{0,1\}$ for every $v\in V$,
\begin{equation*}
|O|=|V|-\sum_{v\in V}w_v=z,\qquad
\sum_{v\in V\setminus O}d(v,C)=\sum_{v\in V}w_vd(v,C).
\end{equation*}
Thus, the minimization over $O$ and the binary linear program over $w$ have the same optimal value, proving the first equality in \eqref{eq:special-cases-outlier-retained-lp}. To verify the second equality in \eqref{eq:special-cases-outlier-retained-lp}, note that
\begin{equation*}
\left\{w\in[0,1]^V:\sum_{v\in V}w_v=|V|-z\right\}
=\conv\left\{w\in\{0,1\}^V:\sum_{v\in V}w_v=|V|-z\right\}.
\end{equation*}
The objective of the linear program in \eqref{eq:special-cases-outlier-retained-lp} is linear, so its minimum is attained at an extreme point in the convex hull above. This proves the second equality in \eqref{eq:special-cases-outlier-retained-lp}.

We next derive the dual of the linear program in \eqref{eq:special-cases-outlier-retained-lp}. Let $t\in\mathbb{R}$ be the multiplier for the equality constraint $\sum_{v\in V}w_v=|V|-z$. The Lagrangian is
\begin{equation*}
\sum_{v\in V}w_vd(v,C)+t\left(|V|-z-\sum_{v\in V}w_v\right)
=(|V|-z)t+\sum_{v\in V}w_v\bigl(d(v,C)-t\bigr).
\end{equation*}
For each $v\in V$, the minimum over $0\leq w_v\leq1$ is attained at $w_v=0$ when $t\leq d(v,C)$ and at $w_v=1$ when $t>d(v,C)$. Hence,
\begin{equation*}
\min_{0\leq w_v\leq1}w_v\bigl(d(v,C)-t\bigr)
=\min\{0,d(v,C)-t\}
=-\bigl(t-d(v,C)\bigr)_{+}.
\end{equation*}
The constraints $0\leq w_v\leq1$ are independent across $v\in V$, so minimizing the Lagrangian over $w_v\in[0,1]$ for $v\in V$ gives
\begin{equation*}
(|V|-z)t+\sum_{v\in V}\min_{0\leq w_v\leq1}w_v\bigl(d(v,C)-t\bigr)
=(|V|-z)t-\sum_{v\in V}\bigl(t-d(v,C)\bigr)_{+}.
\end{equation*}
The linear program in \eqref{eq:special-cases-outlier-retained-lp} is feasible: a binary vector with $|V|-z$ entries equal to one satisfies its constraints. Its objective is nonnegative because $w_v\geq0$ and $d(v,C)\geq0$ for every $v\in V$. Strong duality for linear programming therefore gives
\begin{equation*}
\min_{\substack{w_v\in[0,1]\ (v\in V)\\\sum_{v\in V}w_v=|V|-z}}\sum_{v\in V}w_vd(v,C)
=\max_{t\in\mathbb{R}}\left\{(|V|-z)t-\sum_{v\in V}\bigl(t-d(v,C)\bigr)_{+}\right\}.
\end{equation*}
Combining strong duality with \eqref{eq:special-cases-outlier-retained-lp} proves \eqref{eq:special-cases-outlier-duality}.

It remains to verify that the maximum is attained at $t=d(v,C)$ for some $v\in V$. Let $O^{*}\subseteq V$ be an optimal solution to the minimization on the left-hand side of \eqref{eq:special-cases-outlier-duality}. Such a solution exists because $V$ is finite and $0\leq z<|V|$. Set $t^{*}=\max_{v\in V\setminus O^{*}}d(v,C)$. By the definition of $t^{*}$, we have $d(v,C)\leq t^{*}$ for every $v\in V\setminus O^{*}$. Moreover, $d(v,C)\geq t^{*}$ for every $v\in O^{*}$. Therefore,
\begin{equation}
\label{eq:special-cases-outlier-positive-part}
\bigl(t^{*}-d(v,C)\bigr)_{+}=\begin{cases}
0,&v\in O^{*},\\
t^{*}-d(v,C),&v\in V\setminus O^{*}.
\end{cases}
\end{equation}
It follows that the objective value of the maximization on the right-hand side of \eqref{eq:special-cases-outlier-duality} at $t^{*}$ satisfies
\begin{equation*}
\begin{aligned}
(|V|-z)t^{*}-\sum_{v\in V}\bigl(t^{*}-d(v,C)\bigr)_{+}
&=(|V|-z)t^{*}-\sum_{v\in V\setminus O^{*}}\bigl(t^{*}-d(v,C)\bigr)\\
&=\sum_{v\in V\setminus O^{*}}d(v,C)
=\min_{\substack{O\subseteq V\\|O|=z}}\sum_{v\in V\setminus O}d(v,C),
\end{aligned}
\end{equation*}
where the first equality follows from \eqref{eq:special-cases-outlier-positive-part}, the second equality is due to $|V\setminus O^{*}|=|V|-z$, and the last equality follows from the optimality of $O^{*}$. Thus $t^{*}$ attains the maximum. Since $z<|V|$, the set $V\setminus O^{*}$ is nonempty, and $t^{*}=d(v,C)$ for some $v\in V\setminus O^{*}$. This completes the proof. 
\end{proof}

\begin{proof}[Proof of Proposition~\ref{prop:special-cases-k-median-outliers}]
We first apply Lemma~\ref{lem:special-cases-outlier-duality} to rewrite the minimization over $O$. We then follow an argument similar to the proof of Proposition~\ref{prop:special-cases-robust-k-median} to express the worst-case expectation as a maximum over the target points.

Fix $C\in\Theta$. By Lemma \ref{lem:special-cases-outlier-duality} \eqref{eq:special-cases-outlier-duality} and the definition of $f_C$ in \eqref{eq:special-cases-outlier-threshold-cost}, we have
\begin{equation}
\label{eq:special-cases-outlier-reduced-minimum}
\min_{\substack{O\subseteq V\\|O|=z}}\sum_{v\in V\setminus O}d(v,C)
=\max_{t\in\mathbb{R}}f_C(t).
\end{equation}
Lemma~\ref{lem:special-cases-outlier-duality} also shows that $f_C$ attains its maximum over $\mathbb{R}$ at $t^{*}=d(v,C)$ for some $v\in V$. Since $C$ is a nonempty subset of the finite set $V$, there is a $c\in C$ such that $d(v,C)=d(v,c)$. Hence $t^{*}\in\mathcal{X}$ by the definition of $\mathcal{X}$ in \eqref{eq:special-cases-outlier-threshold-space}. Since $\mathcal{X}\subseteq\mathbb{R}$,
\begin{equation*}
\max_{t\in\mathcal{X}}f_C(t)\leq\max_{t\in\mathbb{R}}f_C(t)=f_C(t^{*})\leq\max_{t\in\mathcal{X}}f_C(t).
\end{equation*}
Combining \eqref{eq:special-cases-outlier-reduced-minimum} with $\max_{t\in\mathcal{X}}f_C(t)=\max_{t\in\mathbb{R}}f_C(t)$ gives
\begin{equation}
\label{eq:special-cases-outlier-discrete-minimum}
    \min_{\substack{O\subseteq V\\|O|=z}}\sum_{v\in V\setminus O}d(v,C)=\max_{t\in\mathcal{X}}f_C(t).
\end{equation}
It remains to express $\max_{t\in\mathcal{X}}f_C(t)$ as the inner worst-case expectation in \eqref{eq:special-cases-gwdro}.

By the definition of the GW ambiguity set in \eqref{eq:finite-gw-ambiguity-set}, the GW ambiguity set for the metric spaces in Proposition~\ref{prop:special-cases-k-median-outliers} is
\begin{equation}
\label{eq:special-cases-outlier-finite-ambiguity-set}
\mathcal{U}_{\rho}^{\mathrm{GW}}(\widehat{\mu}_{\mathcal{Y}})
=\left\{\pi\mathbf{1}_n:\pi\in\mathbb{R}_{+}^{m\times n},\,\pi^{\top}\mathbf{1}_m=\widehat{\mu}_{\mathcal{Y}},\,\sum_{\substack{i,j\in[m]\\k,l\in[n]}}\bigl(|t_i-t_j|-d_{\mathcal{Y}}(y_k,y_l)\bigr)^2\pi_{ik}\pi_{jl}\leq\rho^2\right\}.
\end{equation}
We first show that the GW ambiguity set contains every distribution concentrated on one of the points $t_1,\ldots,t_m$:
\begin{equation*}
\{e_i:i\in[m]\}\subseteq\mathcal{U}_{\rho}^{\mathrm{GW}}(\widehat{\mu}_{\mathcal{Y}}),
\end{equation*}
where $e_i\in\mathbb{R}^{m}$ has entry $1$ at index $i$ and entry $0$ elsewhere, and represents the distribution assigning probability $1$ to $t_i$. Fix $i\in[m]$ and define the coupling $\pi\in\mathbb{R}_{+}^{m\times n}$ as in \eqref{eq:special-cases-point-mass-coupling}, with $G$ replaced by $m$:
\begin{equation}
\label{eq:special-cases-outlier-point-mass-coupling}
\pi_{jl}=\begin{cases}
\widehat{\mu}_{\mathcal{Y},l},&j=i,\\
0,&j\neq i,
\end{cases}
\qquad j\in[m],\quad l\in[n].
\end{equation}
Since $\widehat{\mu}_{\mathcal{Y}}\in\Delta_n$, the coupling satisfies $\pi^{\top}\mathbf{1}_m=\widehat{\mu}_{\mathcal{Y}}$ and $\pi\mathbf{1}_n=e_i$. Substituting \eqref{eq:special-cases-outlier-point-mass-coupling} into the distortion in \eqref{eq:gw-distortion} gives
\begin{equation*}
\mathcal{D}(\pi)
=\sum_{j,j'\in[m]}\sum_{k,l\in[n]}\bigl(|t_j-t_{j'}|-d_{\mathcal{Y}}(y_k,y_l)\bigr)^2\pi_{jk}\pi_{j'l}
=\sum_{k,l\in[n]}\widehat{\mu}_{\mathcal{Y},k}\widehat{\mu}_{\mathcal{Y},l}\bigl(D^{\mathcal{Y}}_{kl}\bigr)^2
=\rho(\widehat{\mu}_{\mathcal{Y}})^2\leq\rho^2.
\end{equation*}
The second equality follows from \eqref{eq:special-cases-outlier-point-mass-coupling}, since $\pi_{ik}=\widehat{\mu}_{\mathcal{Y},k}$ for $k\in[n]$, $\pi_{jk}=0$ for $j\in[m]\setminus\{i\}$ and $k\in[n]$, and $|t_i-t_i|=0$. The third equality follows from the definition of $\rho(\widehat{\mu}_{\mathcal{Y}})$ in \eqref{eq:special-cases-source-radius}, and the inequality follows from the assumption on $\rho$ in Proposition~\ref{prop:special-cases-k-median-outliers}. Thus, by the expression for the GW ambiguity set in \eqref{eq:special-cases-outlier-finite-ambiguity-set}, $e_i=\pi\mathbf{1}_n\in\mathcal{U}_{\rho}^{\mathrm{GW}}(\widehat{\mu}_{\mathcal{Y}})$ for every $i\in[m]$.

We next show that, for the fixed $C\in\Theta$,
\begin{equation}
\label{eq:special-cases-outlier-threshold-identity}
\max_{\mu\in\mathcal{U}_{\rho}^{\mathrm{GW}}(\widehat{\mu}_{\mathcal{Y}})}\sum_{i\in[m]}\mu_i f_C(t_i)
=\max_{t\in\mathcal{X}}f_C(t).
\end{equation}
We first show that the left-hand side of \eqref{eq:special-cases-outlier-threshold-identity} is at least the right-hand side. Since $\mathcal{X}=\{t_i:i\in[m]\}$ and $e_i\in\mathcal{U}_{\rho}^{\mathrm{GW}}(\widehat{\mu}_{\mathcal{Y}})$ for every $i\in[m]$,
\begin{equation*}
\max_{t\in\mathcal{X}}f_C(t)
=\max_{i\in[m]}\sum_{j\in[m]}(e_i)_j f_C(t_j)
\leq\max_{\mu\in\mathcal{U}_{\rho}^{\mathrm{GW}}(\widehat{\mu}_{\mathcal{Y}})}\sum_{j\in[m]}\mu_j f_C(t_j).
\end{equation*}
We next show that the left-hand side of \eqref{eq:special-cases-outlier-threshold-identity} is at most the right-hand side. Every $\mu\in\mathcal{U}_{\rho}^{\mathrm{GW}}(\widehat{\mu}_{\mathcal{Y}})$ satisfies $\mu_j\geq0$ for $j\in[m]$ and $\sum_{j\in[m]}\mu_j=1$. Therefore,
\begin{equation*}
\sum_{j\in[m]}\mu_j f_C(t_j)
\leq\sum_{j\in[m]}\mu_j\max_{t\in\mathcal{X}}f_C(t)
=\max_{t\in\mathcal{X}}f_C(t).
\end{equation*}
Taking the maximum over $\mu\in\mathcal{U}_{\rho}^{\mathrm{GW}}(\widehat{\mu}_{\mathcal{Y}})$ establishes the inequality $\leq$ in \eqref{eq:special-cases-outlier-threshold-identity}. Combining \eqref{eq:special-cases-outlier-threshold-identity} and \eqref{eq:special-cases-outlier-discrete-minimum}, and minimizing over $C\in\Theta$, proves the equivalence of \eqref{eq:special-cases-gwdro} and \eqref{eq:special-cases-k-median-outliers}.
\end{proof}

\subsection{Proofs for Section~\ref{sec:worst-case-upper-bound}}
\label{app:proofs-dnn}

\begin{proof}[Proof of Proposition~\ref{prop:dnn-sdp-bounds}]
If $z$ is feasible for \eqref{eq:worst-case-qcqp}, then setting $Z := zz^{\top}$ gives a feasible solution of \eqref{eq:dnn-relaxation} with the same objective value. Hence, $\vartheta \leq \vartheta_{\mathrm{DNN}}$. Moreover, every feasible solution of \eqref{eq:dnn-relaxation} is feasible for \eqref{eq:sdp-relaxation}. Therefore, $\vartheta_{\mathrm{DNN}} \leq \vartheta_{\mathrm{SDP}}$.
\end{proof}

\begin{proof}[Proof of Corollary~\ref{cor:dnn-exact}]
Under the stated condition, $\bar{z}$ is feasible for \eqref{eq:worst-case-qcqp}. Hence, $\vartheta \geq c^{\top}\bar{z} = \vartheta_{\mathrm{DNN}}$. The reverse inequality follows from Proposition~\ref{prop:dnn-sdp-bounds}.
\end{proof}

\begin{proof}[Proof of Theorem~\ref{thm:facially-reduced-dnn-duality}]
To the optimization instance \eqref{eq:facially-reduced-dnn} we introduce the Lagrangian variables $\lambda\geq0$ associated with the GW constraint, $\tau\in\mathbb{R}$ associated with the equality constraint, and $S \ge 0$ associated with the constraint $VYV^{\top}\geq0$.
The Lagrangian of \eqref{eq:facially-reduced-dnn} is subsequently given by
\begin{equation*}
\lambda\rho^2+\tau
+
\left\langle
\overline{C}
-\lambda\overline{L}
-\tau ee^{\top}
+V^{\top}SV,\,
Y
\right\rangle.
\end{equation*}
The supremum of the above expression over $Y\succeq0$ is finite if and only if
\begin{equation*}
\overline{C}
-\lambda\overline{L}
-\tau ee^{\top}
+V^{\top}SV
\preceq0,
\end{equation*}
in which case the supremum equals $\lambda\rho^2+\tau$. This gives
\eqref{eq:facially-reduced-dnn-dual}.

We then prove strong duality by verifying that Slater's condition holds for the dual \eqref{eq:facially-reduced-dnn-dual}.  Recall $A$ from \eqref{eq:marginal-matrices} and the choice of $F$ and
$z^{0}$ preceding \eqref{eq:facial-reduction-particular-point}. Since
$\mathbf{1}_{n}^{\top}A=\mathbf{1}_{mn}^{\top}$, $AF=0$, and
$Az^{0}=\widehat{\mu}_{\mathcal{Y}}$, we have
\begin{equation*}
\mathbf{1}_{mn}^{\top}F=0,
\qquad
\mathbf{1}_{mn}^{\top}z^{0}
=
\mathbf{1}_{n}^{\top}\widehat{\mu}_{\mathcal{Y}}
=
1.
\end{equation*}
The definition of $V$ in \eqref{eq:facial-reduction-basis} gives $
V^{\top}\mathbf{1}_{mn+1}=2e$.  Moreover, $V$ has full column rank, so $V^{\top}V\succ0$.  Fix $\lambda_{0}>0$ and choose $\eta>0$ sufficiently large that
\begin{equation*}
\lambda_{0}\overline{L}
-\overline{C}
+\eta V^{\top}V
\succ0.
\end{equation*}
Choose $\kappa>\eta$ and set
\begin{equation*}
S_{0}
:=
\kappa\mathbf{1}_{mn+1}\mathbf{1}_{mn+1}^{\top}
-\eta I_{mn+1},
\qquad
\tau_{0}:=4\kappa.
\end{equation*}
Every entry of $S_{0}$ is strictly positive, and
\begin{equation*}
V^{\top}S_{0}V
=
4\kappa ee^{\top}
-\eta V^{\top}V.
\end{equation*}
Consequently,
\begin{equation*}
\lambda_{0}\overline{L}
+\tau_{0}ee^{\top}
-\overline{C}
-V^{\top}S_{0}V =
\lambda_{0}\overline{L}
-\overline{C}
+\eta V^{\top}V
\succ0.
\end{equation*}
Thus, $\lambda_0$, $\tau_0$ and $S_0$ are strictly feasible in \eqref{eq:facially-reduced-dnn-dual}.

Since $(\lambda_0,\tau_0,S_0)$ is feasible for \eqref{eq:facially-reduced-dnn-dual}, weak duality gives
\begin{equation*}
\vartheta_{\mathrm{DNN}}
\leq \lambda_0\rho^2+\tau_0<\infty.
\end{equation*}
The primal problem \eqref{eq:facially-reduced-dnn} is feasible by assumption, so its optimal value is finite.  Since the dual problem is strictly feasible, by Slater's condition we have $\vartheta_{\mathrm{DNN}}=d_{\mathrm{DNN}}$.
\end{proof}

\FloatBarrier

\end{document}